\documentclass{article}
\usepackage{amsmath,amssymb}
\begin{document}
\newtheorem{thm}{Theorem}[section]
\newtheorem{prop}[thm]{Proposition}
\newenvironment{dfn}{\medskip\refstepcounter{thm}
\noindent{\bf Definition \thesection.\arabic{thm}\ }}{\medskip}
\newenvironment{cond}{\medskip\refstepcounter{thm}
\noindent{\bf Condition \thesection.\arabic{thm}\ }}{\medskip}
\newenvironment{ex}{\medskip\refstepcounter{thm}
\noindent{\bf Example \thesection.\arabic{thm}\ }}{\medskip}
\newenvironment{rlist}{\begin{list}{$({\rm \roman{enumi}})$}
{\usecounter{enumi} \setlength{\rightmargin}{10pt}
\setlength{\leftmargin}{40pt} \setlength{\itemsep}{2pt}
\setlength{\parsep}{0pt}
\setlength{\labelwidth}{40pt}}}{\end{list}}
\newenvironment{proof}[1][,]{\medskip\ifcat,#1
\noindent{{\it Proof}:\ }\else\noindent{\it Proof of #1.\ }\fi}
{\hfill$\square$\medskip}
\def\eq#1{{\rm(\ref{#1})}}
\def\P{{\mathbb P}}
\def\Q{{\mathbb Q}}
\def\R{{\mathbb R}}
\def\O{{\mathbb O}}
\def\Z{{\mathbb Z}}
\def\GL{\mathop{\rm GL}}
\def\SL{\mathop{\rm SL}}
\def\SU{\mathop{\rm SU}}
\def\SO{\mathop{\rm SO}}
\def\U{\mathbin{\rm U}}
\def\sech{\mathop{\rm sech}\nolimits}
\def\sn{\mathop{\rm sn}\nolimits}
\def\cn{\mathop{\rm cn}\nolimits}
\def\dn{\mathop{\rm dn}\nolimits}
\def\d{{\rm d}}
\def\w{\wedge}
\def\C{{\mathbb C}}
\def\Gr{\mathop{\rm Gr}\nolimits}
\def\Re{\mathop{\rm Re}\nolimits}
\def\Im{\mathop{\rm Im}\nolimits}
\title{2-Ruled Calibrated 4-folds in $\R^7$ and $\R^8$}
\author{Jason Lotay, University College, Oxford}
\date{}
\maketitle

\linespread{1.1}



\section{Introduction}

In this paper we study certain \emph{calibrated} 4-folds in $\R^7$
and $\R^8$, which are known as \emph{coassociative} in $\R^7$, and
are called \emph{special Lagrangian} (SL) and \emph{Cayley} in
$\R^8$.  We introduce the notion of \emph{2-ruled} 4-folds in
$\R^n$; that is, submanifolds $M$ of $\R^n$ admitting a fibration
$\pi:M\rightarrow\Sigma$ over some 2-fold $\Sigma$ such that each
fibre, $\pi^{-1}(\sigma)$ for $\sigma\in\Sigma$, is an affine
2-plane in $\R^n$.  We say that $M$ is \emph{r-framed} if we are
given an oriented basis for each fibre in a smooth manner.  In such
circumstances there exist orthogonal smooth maps
$\phi_1,\phi_2:\Sigma\rightarrow\mathcal{S}^{n-1}$ and a smooth map
$\psi:\Sigma\rightarrow\R^n$ such that
\begin{equation*}
M =
\{r_1\phi_1(\sigma)+r_2\phi_2(\sigma)+\psi(\sigma):\sigma\in\Sigma,
\hspace{2pt} r_1,r_2\in\R\}.
\end{equation*}
Then the \emph{asymptotic cone} $M_0$ of $M$ is given by:
\begin{equation*}
M_0 = \{r_1\phi_1(\sigma)+r_2\phi_2(\sigma):\sigma\in\Sigma,
\hspace{2pt} r_1,r_2\in\R\}.
\end{equation*}

The motivation for this paper comes from the study of \emph{ruled}
SL 3-folds in $\C^3$ in \cite{Joyce2} and ruled associative
3-folds in $\R^7$ in \cite{Lotay}: these are calibrated 3-folds
that are fibred over a 2-fold by (real) affine straight lines.

We begin in $\S$\ref{calib} by discussing calibrated geometry in
$\R^7$ and $\R^8$.  In particular, we show that coassociative and
SL 4-folds can be considered as special cases of Cayley 4-folds.
In $\S$\ref{2r} we give the definitions required to study 2-ruled
submanifolds.

In $\S$\ref{evolve} we give our main result, Theorem
\ref{Cayevolve}, which is on non-planar, r-framed,\linebreak 2-ruled
Cayley 4-folds. This result characterizes the Cayley condition in
terms of a coupled system of nonlinear first-order partial
differential equations that $\phi_1$ and $\phi_2$ satisfy, and
another such equation on $\psi$ which is \emph{linear} in $\psi$.
Therefore, for a fixed non-planar, r-framed, 2-ruled Cayley cone
$M_0$, the space of r-framed 2-ruled Cayley 4-folds $M$ which have
asymptotic cone $M_0$ has the structure of a finite-dimensional
vector space.

Theorem \ref{Cayholo} gives a means of constructing 2-ruled Cayley
4-folds $M$ from a\linebreak 2-ruled Cayley cone $M_0$, satisfying a
certain condition, involving \emph{holomorphic vector fields}. Using
Theorem \ref{Cayevolve} and Theorem \ref{Cayholo}, we deduce
corresponding results for SL and coassociative 4-folds.

Finally, in $\S$\ref{exs}, we give explicit examples of 2-ruled
4-folds.  We use Theorem \ref{Cayholo} to construct
$\U(1)$-invariant 2-ruled Cayley 4-folds from a
$\U(1)^3$-invariant 2-ruled Cayley cone.  Our other examples are
based on ruled 3-folds and complex cones.

\textit{Acknowledgements.}  I am indebted to my research supervisor,
Dominic Joyce, for his help and guidance.  I would also like to
thank the referees for their hard work in providing a detailed report
and, in particular, for their helpful suggestions regarding the
presentation of the results in Section \ref{ruleds}.

\section{Calibrated 4-folds in $\R^7$ and $\R^8$}
\label{calib}

We begin by defining the basic concepts of \emph{calibrated
geometry} following the approach in \cite{HarLaw}.  Manifolds are
assumed to be smooth and nonsingular everywhere unless stated
otherwise and submanifolds are considered to be immersed.

\begin{dfn}
Let $(M,g)$ be a Riemannian manifold.  An \emph{oriented tangent
$k$-plane} $V$ on $M$ is an oriented $k$-dimensional vector
subspace $V$ of $T_x M$, for some $x$ in $M$. Given an oriented
tangent $k$-plane $V$ on $M$, $g|_V$ is a Euclidean metric on $V$
and hence, using $g|_V$ and the given orientation on $V$, we have
a natural volume form vol$_V$ on $V$ which is a $k$-form on $V$.

 Let $\eta$ be a closed $k$-form on $M$.  Then $\eta$ is a
\emph{calibration} on $M$ if $\eta|_V\leq$ vol$_V$ for all oriented
tangent $k$-planes $V$ on $M$, where
$\eta|_V=\alpha\cdot\text{vol}_V$ for some $\alpha\in\R$, and so
$\eta|_V \leq$ vol$_V$ if $\alpha\leq 1$. Let $N$ be an oriented
$k$-dimensional submanifold of $M$. Then $T_x N$ is an oriented
tangent $k$-plane for all $x\in N$.  We say that $N$ is a
\emph{calibrated submanifold} or \emph{$\eta$-submanifold} if
$\eta|_{T_xN}=$ vol$_{T_x N}$ for all $x\in N$.
\end{dfn}

Calibrated submanifolds are \emph{minimal} submanifolds
\cite[Theorem II.4.2]{HarLaw}.  We define calibrations on $\R^7$
and $\R^8$ as in \cite[Chapter X]{Joyce1}.

\begin{dfn} Let $(x_1,\ldots,x_7)$ be coordinates on $\R^7$ and write
$d{\bf x}_{ij\ldots k}$ for the form $dx_i\w dx_j\w\ldots\w dx_k$.
Define a 3-form $\varphi$ by:
\begin{equation}
\label{phi} \varphi = d{\bf x}_{123}+d{\bf x}_{145}+d{\bf
x}_{167}+d{\bf x}_{246}- d{\bf x}_{257}-d{\bf x}_{347}-d{\bf
x}_{356}.
\end{equation}

\noindent The 4-form $\ast\varphi$, where $\varphi$ and
$\ast\varphi$ are related by the Hodge star, is given by:
\begin{equation}
\label{starphi} \ast\varphi = d{\bf x}_{4567}+d{\bf
x}_{2367}+d{\bf x}_{2345}+d{\bf x}_{1357}-d{\bf x}_{1346}-d{\bf
x}_{1256}-d{\bf x}_{1247}.
\end{equation}

By \cite[Theorem IV.1.16]{HarLaw}, $\ast\varphi$ is a calibration on
$\R^7$ and submanifolds calibrated with respect to $\ast\varphi$ are
called \emph{coassociative} 4-folds.
\end{dfn}

The subgroup of $\GL(7,\R)$ preserving $\varphi$ is $\text{G}_2$.
It is a compact, connected, simply connected, simple,
14-dimensional Lie group, which also preserves the Euclidean
metric on $\R^7$, the orientation on $\R^7$ and $\ast\varphi$. We
note that, by \cite[Theorem IV.1.4]{HarLaw}, $\varphi$ is a
calibration on $\R^7$ and $\varphi$-submanifolds are called
\emph{associative} 3-folds.

\begin{dfn} Let $(x_1,\ldots,x_8)$ be coordinates on $\R^8$.  Define
 a 4-form $\Phi$ by:
\begin{align}
\Phi = &\; d{\bf x}_{1234}+d{\bf x}_{1256}+d{\bf x}_{1278}+d{\bf
x}_{1357}-d{\bf x}_{1368}-d{\bf x}_{1458}-d{\bf x}_{1467}
\nonumber \\
\label{Phi} {}+ &\; d{\bf x}_{5678}+d{\bf x}_{3478}+d{\bf
x}_{3456}+d{\bf x}_{2468}-d{\bf x}_{2457}-d{\bf x}_{2367}-d{\bf
x}_{2358}.
\end{align}

By \cite[Theorem IV.1.24]{HarLaw}, $\Phi$ is a calibration on $\R^8$
and submanifolds calibrated with respect to $\Phi$ are called
\emph{Cayley} 4-folds.
\end{dfn}

The subgroup of $\GL(8,\R)$ preserving $\Phi$ is Spin$(7)$. It is
a compact, connected, simply connected, simple, 21-dimensional Lie
group, which preserves the Euclidean metric and the orientation on
$\R^8$.  It is isomorphic to the double cover of $\SO(7)$.

\medskip

It is worth noting that our formulae \eq{phi}-\eq{Phi} are in
agreement with \cite{Bryant1} and \cite{Joyce1} but differ from
those given in \cite{HarLaw}.  However, they are equivalent up to a
coordinate transformation and a possible change of sign, relating to
a reversal of orientation of the calibrated submanifold, so the
results from \cite{HarLaw} may be applied. Certain formulae used
later from \cite{HarLaw} then need slight modification in order to
be consistent with this choice of representation. These changes
shall be pointed out to the reader.

\medskip

We may consider $\R^8$ as $\C^4$ with complex coordinates
$z_1=x_1+ix_2$, $z_2=x_3+ix_4$, $z_3=x_5+ix_6$, $z_4=x_7+ix_8$, so
we define a calibration on $\C^4$ which can easily be generalised
to $\C^m$.

\begin{dfn} Let $\C^4$ have complex coordinates $(z_1,z_2,z_3,z_4)$ and metric
$g=|dz_1|^2+\ldots+|dz_4|^2$.  Define a real 2-form $\omega$ and a
complex 4-form $\Omega$ on $\C^4$ by:
\begin{equation}
\label{Omega} \omega =  \frac{i}{2}(dz_1\w
d\bar{z}_1+\ldots+dz_4\w d\bar{z}_4),\hspace{10pt} \Omega  =
dz_1\w dz_2\w dz_3\w dz_4.
\end{equation}
Let $L$ be a real oriented 4-fold in $\C^4$.  Then $L$ is a
\emph{special Lagrangian} (SL) 4-fold in $\C^4$ with \emph{phase}
$e^{i\theta}$ if $L$ is calibrated with respect to $\cos\theta\Re
\Omega +\sin\theta \Im\Omega$.  If the phase of $L$ is unspecified
it is taken to be one.
\end{dfn}

Alternative characterizations of these calibrated 4-folds in
$\R^7$ and $\R^8$ are given in \cite{HarLaw}.  The first follows
from \cite[Proposition IV.4.5 \& Theorem IV.4.6]{HarLaw}.

\begin{prop} Let $M$ be a 4-fold in $\R^7$.  Then $M$, with an
appropriate orientation, is coassociative if and only if $\varphi
|_M\equiv 0$.
\end{prop}

The second, for SL 4-folds, is taken from \cite[Corollary
III.1.11]{HarLaw}.

\begin{prop} Let $L$ be a real 4-fold in $\C^4$.  Then $L$, with
the correct orientation, is an SL 4-fold in $\C^4$ with phase
$e^{i\theta}$ if and only if $\omega |_L\equiv 0$ and
$(\sin\theta\Re\Omega-\cos\theta\Im\Omega)|_L\equiv 0$.
\end{prop}

The final result is taken from \cite[Corollary IV.1.29]{HarLaw}.
It requires the definition of the \emph{fourfold cross product} of
four vectors in $\R^8$, for which we identify $\R^8$ with the
\emph{octonions}, or \emph{Cayley numbers}, $\O$.

\begin{dfn}\label{defcross} Let $x,y,z,w\in\O\cong\R^8$.  We define the
\emph{triple cross product} of $x,y,z$ by:
\begin{equation}
\label{triplecross} x\times y\times z =
-\frac{1}{2}\left(x(\bar{y}z)-z(\bar{y}x)\right).
\end{equation}

\noindent The \emph{fourfold cross product} of $x,y,z,w$ is given
by:
\begin{equation*}
\label{fourcross} x\times y\times z\times w =
\frac{1}{4}\big(\bar{x}(y\times z\times w)+ \bar{y}(z\times
x\times w)+\bar{z}(x\times y\times w)+ \bar{w}(y\times x\times
z)\big).
\end{equation*}
\end{dfn}

\vspace{-20pt}

Note that $x\times y\times z\times w$ is an alternating
multi-linear form and that these definitions differ in sign from
those given in \cite{HarLaw} because of the choice of $\Phi$ in
\eq{Phi}.

\begin{prop}
\label{Cay4plane} Let $V$ be a 4-plane in $\R^8$ with basis
$(x,y,z,w)$.  Then $V$, with an appropriate orientation, is Cayley
if and only if $\Im (x\times y\times z\times w)=0$.
\end{prop}

In $\S$\ref{evolve} we shall need the following properties of Cayley
4-folds that relate to \emph{real analyticity}.  The first is a
consequence of the minimality of Cayley 4-folds, as\linebreak
discussed in \cite{HarLaw}, and the second is taken from
\cite[Theorem IV.4.3]{HarLaw}.

\begin{thm}
\label{realanal1}  Suppose that $M$ is a Cayley 4-fold in $\R^8$. Then
$M$ is real analytic.
\end{thm}

\begin{thm}
\label{realanal2} Let $N$ be a real analytic 3-fold in $\R^8$.
Then there exists a unique real analytic Cayley 4-fold in $\R^8$
containing $N$.
\end{thm}

We end the section by giving the following result, which shows
that coassociative and SL 4-folds are special cases of Cayley
4-folds in $\R^8$, the proof of which is immediate from equations
\eq{phi}-\eq{Omega}.

\begin{prop}
\label{specialCay} If we consider $\R^8\cong\R\oplus\R^7$, with
$x_1$ as the coordinate on $\R$ and coordinates on $\R^7$ labelled
as $(x_2,\ldots,x_8)$, then
\begin{equation}
\label{7Phi} \Phi = dx_1\w\varphi + \ast\varphi.
\end{equation}
Hence, if $M\subseteq\R^7\subseteq\R^8$ is a Cayley 4-fold, then
$M$ is a coassociative 4-fold.  Conversely, let $N$ be a
coassociative 4-fold and $M=\{0\}\times
N\subseteq\R\oplus\R^7=\R^8$. Then $M$ is a Cayley 4-fold.

Consider $\R^8\cong\C^4$ with $(z_1,z_2,z_3,z_4)=(x_1+ix_2,
x_3+ix_4,x_5+ix_6,x_7+ix_8)$ as complex coordinates. Then
\begin{equation}
\label{CPhi} \Phi = \frac{1}{2}\,\omega\w\omega + \Re
\hspace{2pt}\Omega.
\end{equation}
Hence, if $M\subseteq\C^4\cong\R^8$ is a Cayley 4-fold such that
$\omega|_M\equiv 0$, then $M$ is an SL 4-fold in $\C^4$. Conversely,
if $L$ is an SL 4-fold in $\C^4\cong\R^8$, then $L$ is a Cayley
4-fold.
\end{prop}

\section{2-Ruled 4-folds in $\R^7$ and $\R^8$}
\label{2r}

We begin by defining \emph{2-ruled} 4-folds in $\R^n$. We take a
\emph{cone} to be a submanifold of $\R^n$ which is invariant under
dilations and nonsingular except possibly at 0.

\begin{dfn}
\label{2Ruled} Let $M$ be a 4-fold in $\R^n$.  A \emph{2-ruling}
of $M$ is a pair $(\Sigma,\pi)$, where $\Sigma$ is a 2-dimensional
manifold and $\pi:M\rightarrow\Sigma$ is a smooth map, such that
$\pi^{-1}(\sigma)$ is an affine 2-plane in $\R^n$ for all
$\sigma\in\Sigma$.  The triple $(M,\Sigma,\pi)$ is a
\emph{2-ruled} 4-fold in $\R^n$.

An \emph{r-framing} for a 2-ruling $(\Sigma,\pi)$ of $M$ is a
choice of oriented orthonormal basis, or frame, for the affine
2-plane $\pi^{-1}(\sigma)$ given by the 2-ruling, for each
$\sigma\in\Sigma$, which varies smoothly with $\sigma$.  Then
$(M,\Sigma,\pi)$ with an r-framing is called \emph{r-framed}.

Let $(M,\Sigma,\pi)$ be an r-framed 2-ruled 4-fold in $\R^n$. For
each $\sigma\in\Sigma$, define $(\phi_1(\sigma),\phi_2(\sigma))$
to be the oriented orthonormal basis for $\pi^{-1}(\sigma)$ given
by the\linebreak r-framing. Then $\phi_1,\phi_2:\Sigma\rightarrow
\mathcal{S}^{n-1}$ are smooth maps.  Define
$\psi:\Sigma\rightarrow\R^n$ such that, for all $\sigma\in\Sigma$,
$\psi(\sigma)$ is the unique vector in $\pi^{-1}(\sigma)$
orthogonal to $\phi_1(\sigma)$ and $\phi_2(\sigma)$.  Then $\psi$
is a smooth map and
\begin{equation}
\label{2ruled} M =
\{r_1\phi_1(\sigma)+r_2\phi_2(\sigma)+\psi(\sigma):\sigma\in\Sigma,
\hspace{2pt} r_1,r_2\in\R\}.
\end{equation}

We define the \emph{asymptotic cone} $M_0$ of a 2-ruled 4-fold $M$
as the set of points in planes $\Pi$ including the origin such
that $\Pi$ is parallel to $\pi^{-1}(\sigma)$ for some
$\sigma\in\Sigma$. If $M$ is r-framed, then
\begin{equation}
\label{2ruledcone} M_0 =
\{r_1\phi_1(\sigma)+r_2\phi_2(\sigma):\sigma\in\Sigma,
\hspace{2pt}r_1,r_2\in\R\}
\end{equation}
and is usually a 4-dimensional cone; that is, whenever the map
$\iota:\Sigma\times\mathcal{S}^1\rightarrow\mathcal{S}^{n-1}$
given by $\iota(\sigma,e^{i\theta})=
\cos\theta\phi_1(\sigma)+\sin\theta\phi_2(\sigma)$ is an
immersion.
\end{dfn}

Let $(M,\Sigma,\pi)$ be a 2-ruled 4-fold in $\R^n$.  Let
$$P=\{({\bf v},\sigma)\in\mathcal{S}^{n-1}\times\Sigma\,:\,\text{${\bf v}$ is
a unit vector parallel to $\pi^{-1}(\sigma)$},\,\sigma\in\Sigma\}$$
and let $\pi_P:P\rightarrow\Sigma$ be given by $\pi_P({\bf
v},\sigma)=\sigma$.  Clearly, $\pi_P:P\rightarrow\Sigma$ is an
$\mathcal{S}^1$ bundle over $\Sigma$. Note that $(M,\Sigma,\pi)$
admits an r-framing if and only if this bundle is trivializable.
Therefore, if $M$ is orientable and $\Sigma$ is non-orientable, e.g.
$\Sigma\cong\mathcal{K}$ where $\mathcal{K}$ is the Klein bottle,
then a 2-ruling $(\Sigma,\pi)$ cannot be r-framed. Moreover, if $M$
is r-framed then $M_0$ is not necessarily 4-dimensional.  For
example, if we take $\Sigma=\R^2$ and define $\phi_1,\phi_2,\psi$ by
$\phi_1(x,y)=(1,0,0,0)$, $\phi_2(x,y)=(0,1,0,0)$ and
$\psi(x,y)=(0,0,x,y)$ for $x,y\in\R$, then $M$, as defined by
\eq{2ruled}, is an r-framed 2-ruled 4-fold since $M=\R^4$, but
$M_0=\R^2$.  We also note that any r-framed 2-ruled 4-fold is
defined by three maps $\phi_1,\phi_2,\psi$ as in \eq{2ruled}. We may
thus construct 2-ruled calibrated 4-folds by formulating
\emph{evolution equations} for $\phi_1,\phi_2,\psi$.

We justify the terminology of asymptotic cone as given in Definition
\ref{2Ruled}. To do this we define the term \emph{asymptotically
conical with order $O(r^{\alpha})$}, where $r$ is the radius
function on $\R^n$, as in \cite[Definition 3.5]{Joyce2}.

\begin{dfn}
\label{asym} Let $M$ be a closed submanifold of $\R^n$ and let
$M_0$ be a closed cone in $\R^n$.  Then $M$ is
\emph{asymptotically conical to $M_0$ with order $O(r^{\alpha})$},
for some $\alpha<1$, if there exist $R>0$, a compact subset $K$
of $M$ and a diffeomorphism $\Psi:M_0\setminus\bar{B}_R\rightarrow
M\setminus K$ such that
\begin{equation}
\label{asymconds} |\nabla^k(\Psi({\bf x})-I({\bf
x}))|=O(r^{\alpha-k}) \hspace{10pt} {\rm for} \hspace{2pt}
k=0,1,2,\ldots \hspace{2pt} {\rm as} \hspace{2pt}
r\rightarrow\infty,
\end{equation}
\noindent where $\bar{B}_R$ is the closed ball of radius $R$ in
$\R^n$ and $I:M_0\setminus \bar{B}_R\rightarrow\R^n$ is the
inclusion map.  Here $|\,.\,|$ is calculated using the cone metric
on $M_0\setminus\bar{B}_R$, and $\nabla$ is a combination of the
Levi-Civita connection derived from the cone metric and the flat
connection on $\R^n$, which acts as partial differentiation.
\end{dfn}

Let $M$ be an r-framed 2-ruled 4-fold with asymptotic cone $M_0$.
Writing $M$ and $M_0$ in the form \eq{2ruled} and \eq{2ruledcone}
respectively, we define, for some $R>0$ and compact $K\subseteq M$,
a diffeomorphism $\Psi:M_0\setminus \bar{B}_R\rightarrow M\setminus
K$ by $\Psi(r_1\phi_1(\sigma)+r_2
\phi_2(\sigma))=r_1\phi_1(\sigma)+r_2\phi_2(\sigma)+\psi(\sigma)$
for all $\sigma\in\Sigma$ and $|r_1|^2+|r_2|^2>R^2$.  Then, if
$\Sigma$ is compact so that $\psi$ is bounded, $\Psi$ satisfies
\eq{asymconds} for $\alpha=0$.  Therefore $M$ is asymptotically
conical to $M_0$ with order $O(1)$.


\section{The Partial Differential Equations}
\label{evolve}

We wish to construct 2-ruled calibrated 4-folds in $\R^7$ and $\R^8$
by solving partial differential equations for maps
$\phi_1,\phi_2,\psi$. By Proposition \ref{specialCay}, it is
sufficient to consider the Cayley case.

Let $\Sigma$ be a 2-dimensional, connected, real analytic
manifold, let $\phi_1,\phi_2:\Sigma\rightarrow\mathcal{S}^7$ be
orthogonal real analytic maps such that
$\iota:\Sigma\times\mathcal{S}^1\rightarrow\mathcal{S}^7$ defined
by
$\iota(\sigma,e^{i\theta})=\cos\theta\phi_1(\sigma)+\sin\theta\phi_2(\sigma)$
is an immersion and let $\psi:\Sigma\rightarrow\R^8$ be a real
analytic map.  Clearly, $\R^2\times\Sigma$ is an r-framed 2-ruled
4-fold with 2-ruling $(\Sigma,\pi)$, where
$\pi(r_1,r_2,\sigma)=\sigma$.  Let $M$ be defined by 
\eq{2ruled}.  Then $M$ is the image of the map
$\iota_M:\R^2\times\Sigma\rightarrow\R^8$ given by
$\iota_M(r_1,r_2,\sigma)=r_1\phi_1(\sigma)+r_2\phi_2(\sigma)+
\psi(\sigma)$. Since $\iota$ is an immersion, $\iota_M$ is an
immersion almost everywhere.  Thus $M$ is an r-framed 2-ruled
4-fold in $\R^8$, possibly with singularities.

Suppose that $M$ is Cayley and let $p\in M$.  Then there exist
$(r_1,r_2)\in\R^2$, $\sigma\in\Sigma$ such that
$p=r_1\phi_1(\sigma)+r_2\phi_2(\sigma)+\psi(\sigma)$.  Choose
 coordinates $(s,t)$ near $\sigma$ in $\Sigma$. Then
$T_pM=\langle x,y,z,w\rangle_{\R}$, where $x=\phi_1(\sigma)$,
$y=\phi_2(\sigma)$, $z=r_1\frac{\partial\phi_1}{\partial
s}(\sigma)+ r_2\frac{\partial\phi_2}{\partial
s}(\sigma)+\frac{\partial\psi}{\partial s} (\sigma)$,
$w=r_1\frac{\partial\phi_1}{\partial t}(\sigma)+
r_2\frac{\partial\phi_2}{\partial
t}(\sigma)+\frac{\partial\psi}{\partial t} (\sigma)$. The tangent
space $T_pM$ is a Cayley 4-plane.  By Proposition \ref{Cay4plane},
this is true if and only if $\Im(x\times y\times z\times w)=0$.
This implies that a quadratic in $r_1,r_2$ must vanish, but since
this condition is forced to hold for all $(r_1,r_2)\in \R^2$, each
coefficient in the quadratic is zero.  Therefore the following set
of equations must hold in $\Sigma$:
\begin{align}
\label{Cay1} \Im\left(\phi_1\times\phi_2\times
\frac{\partial\phi_1}{\partial
s}\times\frac{\partial\phi_1}{\partial t}\right)
& =  0, \\
\label{Cay2} \Im\left(\phi_1\times\phi_2\times
\frac{\partial\phi_2}{\partial
s}\times\frac{\partial\phi_2}{\partial t}\right)
& =  0, \\
\label{Cay3} \Im\left(\phi_1\times\phi_2\times
\frac{\partial\phi_1}{\partial
s}\times\frac{\partial\phi_2}{\partial t}\right) +
\Im\left(\phi_1\times\phi_2\times \frac{\partial\phi_2}{\partial
s}\times\frac{\partial\phi_1}{\partial t}\right) & =  0, \\
\label{Cay4} \Im\left(\phi_1\times\phi_2\times
\frac{\partial\psi}{\partial s}\times\frac{\partial\psi}{\partial
t}\right)
& = 0, \\
\label{Cay5} \Im\left(\phi_1\times\phi_2\times
\frac{\partial\phi_1}{\partial
s}\times\frac{\partial\psi}{\partial t}\right) +
\Im\left(\phi_1\times\phi_2\times \frac{\partial\psi}{\partial
s}\times\frac{\partial\phi_1}{\partial t}\right)
& =  0, \\
\label{Cay6} \Im\left(\phi_1\times\phi_2\times
\frac{\partial\phi_2}{\partial
s}\times\frac{\partial\psi}{\partial t}\right)
+\Im\left(\phi_1\times\phi_2\times \frac{\partial\psi}{\partial
s}\times\frac{\partial\phi_2}{\partial t}\right) & =  0.
\end{align}

If we do not suppose $M$ to be Cayley but instead insist that
\eq{Cay1}-\eq{Cay6} hold in $\Sigma$ then, following the argument
above, each tangent space to $M$ must be Cayley and hence $M$ is a
Cayley 4-fold.  Noting that \eq{Cay1}-\eq{Cay3} are precisely the
conditions for the asymptotic cone $M_0$ of $M$ to be Cayley, we
deduce the following result.

\begin{prop}
\label{calibcone} Let $M$ be a 2-ruled Cayley 4-fold in $\R^8$ and
let $M_0$ be its asymptotic cone.  Then $M_0$ is a 2-ruled Cayley
cone in $\R^8$ provided that it is\linebreak 4-dimensional.
\end{prop}

Clearly, $M_0$ is the image of the map
$\iota_0:\R^2\times\Sigma\rightarrow\R^8$ given by
$\iota_0(r_1,r_2,\sigma)=r_1\phi_1(\sigma)+r_2\phi_2(\sigma)$. Since
we suppose that $\iota$ is an immersion, $\iota_0$ is an immersion
except at $(r_1,r_2)=(0,0)$, so $M_0$ is nonsingular except at 0 and
thus is a cone.

Note that $\Phi$ is a nowhere vanishing 4-form on $M_0$ that defines
its orientation, since $M_0$ is Cayley.  Hence, if $(s,t)$ are local
coordinates on $\Sigma$, we can define them to be oriented by
imposing
\begin{equation}
\label{orient1}
\Phi\left(\phi_1,\hspace{2pt}\phi_2,\hspace{2pt}\hspace{2pt}r_1
\frac{\partial\phi_1}{\partial s}+
r_2\frac{\partial\phi_2}{\partial
s},\hspace{2pt}r_1\frac{\partial\phi_1}{\partial
t}+r_2\frac{\partial\phi_2}{\partial t}\right)>0
\end{equation}
\noindent for all $(r_1,r_2)\in\R^2\setminus\{(0,0)\}$.  It
follows that
\begin{equation}
\label{orient2}
\Phi\left(\phi_1,\phi_2,\frac{\partial\phi_1}{\partial
s},\frac{\partial\phi_1}{\partial t}\right)>0 \quad \text{and}
\quad\Phi\left(\phi_1,\phi_2,\frac{\partial\phi_2}{\partial
s},\frac{\partial\phi_2}{\partial t}\right)>0.
\end{equation}
\noindent Consequently,
$\{\phi_1,\phi_2,\frac{\partial\phi_j}{\partial
s},\frac{\partial\phi_j}{\partial t}\}$ is a linearly independent
set for $j=1,2$.  Moreover, \eq{orient1} is equivalent to the
condition that $\iota$ is an immersion.

\medskip

We now construct a metric on $\Sigma$, under suitable conditions,
using $\phi_1,\phi_2$ and the metric on $\R^8$. This enables us to
formulate \eq{Cay1}-\eq{Cay6} as partial differential equations
involving the triple cross product in $\R^8$, which we may write as
\begin{equation}
\label{Caytriplecross} (x\times y\times z)^e =
\Phi_{abcd}x^ay^bz^cg^{de}
\end{equation}
\noindent using index notation for tensors on $\R^8$, where $g^{de}$
is the inverse of the Euclidean metric on $\R^8$ and $x,y,z\in\R^8$.
It can be easily verified, using \eq{Phi} and a
multiplication table for the octonions as given in \ref{octotable},
that this definition coincides with \eq{triplecross}. We immediately
deduce that
\begin{equation}
\label{Caytriplecross2} \Phi(x,y,z,w)=g(x\times y\times z,w).
\end{equation}
\noindent Note that the triple cross product $x\times y\times z$ is
orthogonal to $x,y,z$, and that it is nonzero if and only if
$\{x,y,z\}$ is a linearly independent set.

\medskip

For a function $f:\Sigma\rightarrow\R^8$, we define
$f^{\perp}:\Sigma\rightarrow\R^8$ by choosing $f^{\perp}(\sigma)$ to
be the component of $f(\sigma)$ that lies in the orthogonal
complement of $\langle\phi_1(\sigma),\phi_2(\sigma)\rangle_{\R}$.
Since the fourfold cross product is alternating, \eq{Cay1}-\eq{Cay3}
hold if and only if
\begin{equation}\label{fourfoldtheta}
\Im\left(\phi_1\!\times\phi_2\!\times\!\left(\cos\theta\frac{\partial\phi_1}{\partial
s}^{\perp}\!\!\!+\sin\theta\frac{\partial\phi_2}{\partial
s}^{\perp}\right)\!\times\left(\cos\theta\frac{\partial\phi_1}{\partial
t}^{\perp}\!\!\!+\sin\theta\frac{\partial\phi_2}{\partial
t}^{\perp}\right)\right)=0,
\end{equation}
\noindent for all $\theta\in\R$.  Let $\sigma\in\Sigma$.  From
Proposition \ref{Cay4plane} and \eq{orient1} we see that, for each
$\theta\in\R$, the four terms in \eq{fourfoldtheta}, evaluated at
$\sigma$, form a basis for a Cayley 4-plane $\Pi_{\theta}$. By the
definition of the triple cross product, we may also take
$(\phi_1(\sigma),\phi_2(\sigma),\cos\theta
\frac{\partial\phi_1}{\partial s}^{\perp}(\sigma)
+\sin\theta\frac{\partial\phi_2}{\partial s}^{\perp}(\sigma)
,\cos\theta\phi_1\times\phi_2\times \frac{\partial\phi_1}{\partial
s}^{\perp}(\sigma)
+\sin\theta\phi_1\times\phi_2\times\frac{\partial\phi_2}{\partial
s}^{\perp}(\sigma))$ as a basis for $\Pi_{\theta}$.  Therefore,
\begin{gather}
\cos\theta\frac{\partial\phi_1}{\partial
t}^{\perp}(\sigma)+\sin\theta\frac{\partial\phi_2}{\partial
t}^{\perp}(\sigma)=A_{\theta}\left(\cos\theta\frac{\partial\phi_1}{\partial
s}^{\perp}(\sigma)+\sin\theta\frac{\partial\phi_2}{\partial
s}^{\perp}(\sigma)\right) \nonumber\\
\label{li} {}+B_{\theta}\left(\cos\theta\phi_1\times\phi_2\times
\frac{\partial\phi_1}{\partial s}^{\perp}(\sigma)
+\sin\theta\phi_1\times\phi_2\times\frac{\partial\phi_2}{\partial
s}^{\perp}(\sigma)\right)
\end{gather}
for constants $A_{\theta},B_{\theta}$ depending on $\theta$. We set
$\theta=0,\frac{\pi}{2}$ in \eq{li} and substitute back in the
expressions found for the $t$ derivatives to obtain:
\begin{gather}
\cos\theta\left((A_{0}-A_{\theta})\frac{\partial\phi_1}{\partial
s}^{\perp}(\sigma)+(B_{0}-B_{\theta})\phi_1\times\phi_2\times
\frac{\partial\phi_1}{\partial s}^{\perp}(\sigma)\right)
=\hspace{24pt} \nonumber\\[2pt]
\label{li2}
\hspace{24pt}\sin\theta\left((A_{\theta}-A_{\frac{\pi}{2}})\frac{\partial\phi_2}{\partial
s}^{\perp}(\sigma)
+(B_{\theta}-B_{\frac{\pi}{2}})\phi_1\times\phi_2\times\frac{\partial\phi_2}{\partial
s}^{\perp}(\sigma)\right)
\end{gather}
To proceed in defining a metric on $\Sigma$ we impose a condition on
the dimension of $$V_{\sigma}=\left\langle
\frac{\partial\phi_1}{\partial s}^{\perp}(\sigma),
\frac{\partial\phi_2}{\partial s}^{\perp}(\sigma)
,\phi_1\times\phi_2\times\frac{\partial\phi_1}{\partial
s}^{\perp}(\sigma),\phi_1\times\phi_2\times\frac{\partial\phi_2}{\partial
s}^{\perp}(\sigma)\right\rangle_{\R}.$$

Let $W_{\sigma}=\langle
\phi_1(\sigma),\phi_2(\sigma)\rangle_{\R}^{\perp}\subseteq\R^8$ and
define $J:W_{\sigma}\rightarrow W_{\sigma}$ by
$J(v)=\phi_1(\sigma)\times\phi_2(\sigma)\times v$.  It is clear,
through calculation in coordinates, that $J^2=-1$ on $W_{\sigma}$.
Note that $V_{\sigma}\subseteq W_{\sigma}$ is closed under the
action of $J$, which can thus be considered as a form of complex
structure on $V_{\sigma}$.  Hence, $V_{\sigma}$ is even dimensional.
Since the case dim$\,V_{\sigma}=0$ is excluded by \eq{orient2},
dim$\,V_{\sigma}=2$ or $4$.  Recall that $\Sigma$ is real analytic
and connected.  Therefore $\{\sigma\in\Sigma:\text{dim}\,
V_{\sigma}=2\}$ is a closed real analytic subset of $\Sigma$ and
consequently either coincides with $\Sigma$ or is of zero measure in
$\Sigma$.

Suppose that dim$\,V_{\sigma}=4$.  The four vectors in \eq{li2} are
then linearly independent and hence
\begin{equation*}
(A_0-A_{\theta})\cos\theta=(B_0-B_{\theta})\cos\theta=
(A_{\frac{\pi}{2}}-A_{\theta})\sin\theta=(B_{\frac{\pi}{2}}-B_{\theta})\sin\theta=0
\end{equation*}
\noindent for all $\theta$.  This clearly forces
$A_{\theta},B_{\theta}$ to be constant. Let $A_{\theta}=A$ and
$B_{\theta}=B$ for all $\theta$, where $A,B$ are real constants. We
define a metric $g_{\Sigma}$ on $\Sigma$ pointwise by the following
equations:
\begin{equation}
\label{metric1} g_{\Sigma}\!\left(\frac{\partial}{\partial s}\,,
\frac{\partial}{\partial
t}\right)=Ag_{\Sigma}\!\left(\frac{\partial}{\partial s}\,,
\frac{\partial}{\partial s}\right),\;
g_{\Sigma}\!\left(\frac{\partial}{\partial t}\,,
\frac{\partial}{\partial
t}\right)=(A^2+B^2)g_{\Sigma}\!\left(\frac{\partial}{\partial s}\,,
\frac{\partial}{\partial s}\right).
\end{equation}
Using \eq{li} and the fact that $J^2=-1$ on $V_{\sigma}$,
\begin{equation*}
\left(\begin{array}{c}\phi_1\times\phi_2\times
\frac{\partial\phi_j}{\partial s}^{\perp}(\sigma) \\
\phi_1\times\phi_2\times \frac{\partial\phi_j}{\partial
t}^{\perp}(\sigma)\end{array}\right) = K\left(\begin{array}{c}
\frac{\partial\phi_j}{\partial s}^{\perp}(\sigma) \\
\frac{\partial\phi_j}{\partial
t}^{\perp}(\sigma)\end{array}\right),
\end{equation*}
for $j=1,2$, where $K$ is a $2\times 2$ matrix given by:
\begin{equation*}
K=\frac{1}{B}\left(\begin{array}{cc}-A & 1 \\ -(A^2+B^2) &
A\end{array}\right).
\end{equation*}
If we change coordinates $(s,t)$ to $(\tilde{s},\tilde{t})$, with
Jacobian matrix $L$, then $K$ transforms to a matrix
$\tilde{K}=LKL^{-1}$.  We may then calculate the corresponding
$\tilde{A},\tilde{B}$ defining $\tilde{K}$ and see that they satisfy
\eq{metric1} for the coordinates $(\tilde{s},\tilde{t})$. Therefore,
$g_{\Sigma}$ is a well-defined metric, up to scale, covariant under transformation
of coordinates.

Having defined the metric $g_{\Sigma}$ we can consider $\Sigma$ as a
\emph{Riemannian} 2-fold, which has a natural orientation derived
from the orientation on $M$ and on the 2-planes
$\langle\phi_1(\sigma),\phi_2(\sigma)\rangle_{\R}$. Therefore it has
a natural \emph{complex structure} which we denote as $J$. If we
choose a local holomorphic coordinate $u=s+it$ on $\Sigma$, then the
corresponding real coordinates must satisfy
$\frac{\partial}{\partial t}=J\frac{\partial}{\partial s}$.  We say
that local real coordinates $(s,t)$ on $\Sigma$ satisfying this
condition are \emph{oriented conformal coordinates}.  This forces
$A=0$, $B=1$ in the notation of \eq{metric1}, since $B>0$ by
\eq{orient2}.

\medskip

We now state and prove a theorem in this case.

\begin{thm}
\label{evolve1} Let $\Sigma$ be a connected real analytic 2-fold,
let $\phi_1,\phi_2:\Sigma\rightarrow\mathcal{S}^7$ be orthogonal
real analytic maps such that the map
$\iota:\Sigma\times\mathcal{S}^1\rightarrow\mathcal{S}^7$ defined by
$\iota(\sigma,e^{i\theta})=\cos\theta\phi_1(\sigma)+\sin\theta\phi_2(\sigma)$
is an immersion, and let $\psi:\Sigma\rightarrow\R^8$ be a real
analytic map. Define $M$ by \eq{2ruled} and suppose that
dim$\,V_{\sigma}=4$ almost everywhere in $\Sigma$. Then $M$ is
Cayley if and only if
\begin{align}
\label{Cay7}\frac{\partial\phi_1}{\partial t} & =
\phi_1\times\phi_2\times \frac{\partial\phi_1}{\partial
s} + f\phi_2, \\
\label{Cay8} \frac{\partial\phi_2}{\partial t} & =
\phi_1\times\phi_2\times \frac{\partial\phi_2}{\partial s} - f
\phi_1,
\end{align}
\noindent for some function $f:\Sigma\rightarrow\R$, and $\psi$
satisfies
\begin{equation}
\label{Caypsi} \frac{\partial\psi}{\partial t} =
\phi_1\times\phi_2\times\frac{\partial\psi}{\partial s} +
g_1\phi_1+g_2\phi_2
\end{equation}
\noindent for some functions $g_1,g_2:\Sigma\rightarrow\R$, where
the triple cross product is defined in \eq{Caytriplecross} and
$(s,t)$ are oriented conformal coordinates on $\Sigma$.  Moreover,
sufficiency holds irrespective of the dimension of $V_{\sigma}$.
\end{thm}

\begin{proof} Recalling that \eq{Cay1}-\eq{Cay6} correspond
to the condition that $M$ is Cayley, we show that
\eq{Cay1}-\eq{Cay3} are equivalent to \eq{Cay7}-\eq{Cay8}, and that
\eq{Cay4}-\eq{Cay6} are equivalent to \eq{Caypsi}.

Let $\sigma\in\Sigma$.  Since $\phi_1$ maps to $\mathcal{S}^7$ it is
clear that $\phi_1(\sigma)$ is orthogonal to
$\frac{\partial\phi_1}{\partial s}(\sigma)$ and
$\frac{\partial\phi_1}{\partial t}(\sigma)$.  By \eq{li} and the
work above, there exist $a_1,a_2,a_3\in\R$ such that
\begin{equation}
\label{dt1} \frac{\partial\phi_1}{\partial t}(\sigma) =
a_1\phi_1\times\phi_2\times\frac{\partial\phi_1}{\partial
s}(\sigma)+a_2\phi_2(\sigma)+a_3\frac{\partial\phi_1}{\partial
s}(\sigma).
\end{equation}
We then calculate:
\begin{equation*}
g\left(\frac{\partial\phi_1}{\partial
t}^{\perp}(\sigma),\frac{\partial\phi_1}{\partial
s}^{\perp}(\sigma)\right)=a_3\left|\frac{\partial\phi_1}{\partial
s}^{\perp}(\sigma)\right|^2.
\end{equation*}
\noindent The left-hand side is zero by \eq{metric1} since $(s,t)$
are oriented conformal coordinates, and hence $a_3=0$.  We also have
that
\begin{equation*}
\left|\frac{\partial\phi_1}{\partial
t}^{\perp}(\sigma)\right|^2=a_1^2
\left|\phi_1\times\phi_2\times\frac{\partial\phi_1}{\partial
s}(\sigma)\right|^2= a_1^2\left|\frac{\partial\phi_1}{\partial
s}^{\perp}(\sigma)\right|^2.
\end{equation*}
Therefore $a_1^2=1$ by \eq{metric1}. Further, taking the inner
product of \eq{dt1} with the triple cross product gives:
\begin{align*}
\left|\phi_1\times\phi_2\times\frac{\partial\phi_1}{\partial
s}(\sigma)\right|^2a_1 & =  g\left(\frac{\partial\phi_1}{\partial
t}(\sigma),\phi_1\times\phi_2\times\frac{\partial\phi_1}{\partial
s}(\sigma)\right) \\
& =
\Phi\left(\phi_1(\sigma),\phi_2(\sigma),\frac{\partial\phi_1}{\partial
s}(\sigma),\frac{\partial\phi_1}{\partial t}(\sigma)\right),
\end{align*}
\noindent using \eq{Caytriplecross}. Therefore $a_1>0$ by
\eq{orient2}.  Hence, $a_1=1$ and \eq{Cay7} holds at
$\sigma$ with $f(\sigma)=a_2$.

If \eq{Cay7} holds at $\sigma$ then, by the definition of the triple
cross product, the 4-plane spanned by
$\{\phi_1(\sigma),\phi_2(\sigma),\frac{\partial\phi_1}{\partial
s}(\sigma), \frac{\partial\phi_1}{\partial t}(\sigma)\}$ is Cayley.

Similarly, we deduce that \eq{Cay2} holding at $\sigma$ is
equivalent to
\begin{equation}
\label{phi2} \frac{\partial\phi_2}{\partial t} =
\phi_1\times\phi_2\times\frac{\partial \phi_2}{\partial
s}+f^{\prime}\phi_1
\end{equation}
\noindent at $\sigma$, for some function
$f^{\prime}:\Sigma\rightarrow\R$.  However,
\begin{equation*}
\frac{\partial}{\partial
t}\,g(\phi_1,\phi_2)=g\left(\frac{\partial\phi_1}{\partial t}
,\phi_2\right)+g\left(\phi_1,\frac{\partial\phi_2}{\partial
t}\right)=0
\end{equation*}
\noindent and hence $f^{\prime}=-f$.

It follows from \cite[Theorem IV.1.38]{HarLaw} that Spin$(7)$ acts
transitively upon oriented orthonormal bases of Cayley 4-planes.
Therefore, we are always able to transform coordinates on $\R^8$
using Spin$(7)$ such that a Cayley 4-plane has basis
$(e_1,e_2,e_3,e_4)$. In particular, any orthonormal pair can be
mapped to the pair $(e_1,e_2)$.

By the remarks above, we may transform coordinates on $\R^8$ using
Spin$(7)$ such that $\phi_1(\sigma)=e_1$, $\phi_2(\sigma)=e_2$,
$\frac{\partial\phi_1}{\partial s}(\sigma)=b_1e_1+\ldots+b_8e_8$,
$\frac{\partial\phi_2}{\partial
s}(\sigma)=b_1^{\prime}e_1+\ldots+b_8^{\prime}e_8$ for some real
constants $b_1,\ldots,b_8,b_1^{\prime},\ldots,b_8^{\prime}$. If
\eq{Cay7} and \eq{Cay8} hold, we may calculate
$\frac{\partial\phi_1}{\partial t}(\sigma)$ and
$\frac{\partial\phi_2}{\partial t}(\sigma)$. A straightforward
calculation in coordinates then shows that \eq{Cay1}-\eq{Cay3} hold
at $\sigma$.  Since the triple cross product is invariant under
Spin$(7)$ by equation \eq{Caytriplecross}, we conclude that
\eq{Cay1}-\eq{Cay3} are equivalent to \eq{Cay7} and \eq{Cay8}.

Suppose now that \eq{Cay4}-\eq{Cay6} hold at $\sigma\in\Sigma$.
Using Spin$(7)$, transform coordinates such that
$\phi_1(\sigma)=e_1$, $\phi_2(\sigma)=e_2$,
$\frac{\partial\phi_1}{\partial s}(\sigma)=b_1e_1+\ldots+b_4e_4$,
where $b_1,\ldots,b_4$ are real constants, which we are free to do
by equation \eq{Cay1}.  In these coordinates write
$\frac{\partial\psi}{\partial s}(\sigma)=c_1e_1+\ldots+c_8e_8$ and
$\frac{\partial\psi}{\partial t}(\sigma)=d_1e_1+\ldots+d_8e_8$.
Calculating $\frac{\partial\phi_1}{\partial t}(\sigma)$ using
\eq{Cay7}, we then evaluate the terms in \eq{Cay5} as follows:
\begin{align}
\label{Cay52}\left(\begin{array}{rr} -b_4 & -b_3 \\ -b_3 & b_4
\end{array}\right)\left(\begin{array}{c} d_5+c_6 \\ d_6-c_5
\end{array}\right)
&=0, \\
\label{Cay53}\left(\begin{array}{rr} b_4 & b_3 \\ b_3 & -b_4
\end{array}\right)\left(\begin{array}{c} d_7+c_8 \\
d_8-c_7 \end{array}\right)&=0.
\end{align}
Details of the calculation of the fourfold cross product may be
found in \ref{calc}. The determinant of the matrices in \eq{Cay52}
and \eq{Cay53} is $-b_3^2-b_4^2\neq 0$, since
$\frac{\partial\phi_1}{\partial
s}(\sigma)\notin\langle\phi_1(\sigma),\phi_2(\sigma)\rangle_{\R}$.
Therefore $d_5=-c_6$, $d_6=c_5$, $d_7=-c_8$ and $d_8=c_7$.

We may also evaluate \eq{Cay4}:
\begin{align}
\label{Cay42} c_5d_8+c_6d_7-c_7d_6-c_8d_5&=0, \\
\label{Cay43} c_5d_7-c_6d_8-c_7d_5+c_8d_6&=0, \\
\label{Cay44} -c_3d_8+c_7d_4-c_4d_7+c_8d_3&=0, \\
\label{Cay45} -c_3d_7+c_4d_8-c_8d_4+c_7d_3&=0, \\
\label{Cay46} c_3d_6+c_4d_5-c_5d_4-c_6d_3&=0, \\
\label{Cay47} c_3d_5-c_4d_6-c_5d_3+c_6d_4&=0.
\end{align}
Again, details may be found in \ref{calc}.
 Substituting in the results above, we have that
\eq{Cay42}-\eq{Cay43} are satisfied trivially and
\eq{Cay44}-\eq{Cay47} become:
\begin{align}
\label{Cay48} \left(\begin{array}{rr} c_8 & c_7 \\ c_7 & -c_8
\end{array}\right)\left(\begin{array}{c} d_3+c_4 \\ d_4-c_3
\end{array}\right)
&=0, \\
\label{Cay49} \left(\begin{array}{rr} -c_6 & -c_5 \\ -c_5 & c_6
\end{array}\right)\left(\begin{array}{c} d_3+c_4 \\ d_4-c_3
\end{array}\right)
&=0.
\end{align}
We deduce that the determinants of the matrices in \eq{Cay48} and
\eq{Cay49} are zero, or the vector appearing in both equations is
zero.  Therefore,

\smallskip

\begin{rlist}
\item $d_3=-c_4$ and $d_4=c_3$

\smallskip

\noindent or

\smallskip

\item $c_5=c_6=c_7=c_8=0$.
\end{rlist}
\smallskip

Condition (i) implies that \eq{Caypsi} holds at $\sigma$ with
$g_1(\sigma)=d_1$, $g_2(\sigma)=d_2$, by the definition of the
triple cross product and its invariance under Spin$(7)$.  Condition
(ii) corresponds to
\begin{equation}
\label{Caypsi2} \frac{\partial\psi}{\partial
s}(\sigma),\frac{\partial\psi}{\partial
t}(\sigma)\in\bigg\langle\phi_1(\sigma),\phi_2(\sigma),\frac{\partial\phi_j}{\partial
s}(\sigma),\frac{\partial\phi_j}{\partial
t}(\sigma)\bigg\rangle_{\R}
\end{equation}
holding for $j=1$.  Therefore, \eq{Cay4} and \eq{Cay5} are
equivalent to \eq{Caypsi} or \eq{Caypsi2} for $j=1$ holding at
$\sigma$.  We similarly deduce that \eq{Cay4} and \eq{Cay6} are
equivalent to \eq{Caypsi} or \eq{Caypsi2} for $j=2$ holding at
$\sigma$.

We conclude that \eq{Cay1}-\eq{Cay6} are equivalent to \eq{Cay7},
\eq{Cay8} and condition \eq{Caypsi} or \eq{Caypsi2} for $j=1,2$ at
each point $\sigma\in\Sigma$.
Recall that $\Sigma$ is connected and $\phi_1,\phi_2,\psi$ and
$\Sigma$ are real analytic.  Note that
$\Sigma_1=\{\sigma\in\Sigma:\text{dim}\,V_{\sigma}=4\}$ is an open
subset of $\Sigma$ whose complement is measure zero in $\Sigma$ by
hypothesis.

Let $\sigma\in\Sigma_1$ and suppose that \eq{Caypsi2} for $j=1,2$
holds at $\sigma$. Then, there exist real constants $C_{jk}$, for
$j=1,2$, $1\leq k\leq4$, such that
$$\frac{\partial\psi}{\partial s}(\sigma)=C_{j1}\phi_1(\sigma)+C_{j2}
\phi_2(\sigma)+C_{j3}\frac{\partial\phi_j}{\partial
s}^{\perp}(\sigma)+C_{j4}\frac{\partial\phi_j}{\partial
t}^{\perp}(\sigma).
$$
Clearly, $C_{1k}=C_{2k}$ for $k=1,2$ by the definition of
$g^{\perp}$ for a function $g$.  Moreover, since dim$\,V_{\sigma}=4$
ensures the linear independence of the partial derivatives of
$\phi_1$ and $\phi_2$, $C_{jk}=0$ for $j=1,2$, $k=3,4$.  Hence,
$\frac{\partial\psi}{\partial s}(\sigma)$ and, similarly,
$\frac{\partial\psi}{\partial t}(\sigma)$ lie in
$\langle\phi_1(\sigma),\phi_2(\sigma)\rangle_{\R}$ for almost all
$\sigma\in\Sigma$.  Therefore $\psi$ satisfies \eq{Caypsi}.

Consequently, \eq{Caypsi} holds in $\Sigma_1$.  Moreover,
$\Sigma_2=\{\sigma\in\Sigma:\text{\eq{Caypsi} holds at $\sigma$}\}$
is a closed real analytic subset of $\Sigma$ and so must either
coincide with $\Sigma$ or be of zero measure in $\Sigma$.  Since
$\Sigma_1\subseteq\Sigma_2$, $\Sigma_2$ cannot be measure zero and
so must equal $\Sigma$. This completes the proof.
\end{proof}


Note that \eq{Caypsi} is a \emph{linear} condition on $\psi$ given
$\phi_1$ and $\phi_2$, and that \eq{Cay7} and \eq{Cay8} are
equivalent to the fact that the asymptotic cone $M_0$ of $M$ is
Cayley.  Therefore, if we are given an r-framed, 2-ruled Cayley cone
$M_0$ defined by $\phi_1$ and $\phi_2$, then any solution $\psi$ of
\eq{Caypsi}, together with $\phi_1$ and $\phi_2$, defines an
r-framed 2-ruled Cayley 4-fold with asymptotic cone $M_0$.  We also
note that \eq{Caypsi} is unchanged if $\phi_1$ and $\phi_2$ are
fixed and satisfy \eq{Cay7} and \eq{Cay8}, but $\psi$ is replaced by
$\psi+\tilde{g}_1\phi_1+\tilde{g}_2\phi_2$ for real analytic maps
$\tilde{g}_1,\tilde{g}_2$.  We can thus locally transform $\psi$
such that $g_1$ and $g_2$ are zero.

\medskip

If we suppose instead that dim$\,V_{\sigma}=2$ for all
$\sigma\in\Sigma$ then we are unable, in general, to define a
suitable metric and hence oriented conformal coordinates on
$\Sigma$. However, we shall show that if we exclude \emph{planar}
r-framed 2-ruled 4-folds, then \eq{Cay7}-\eq{Caypsi} of Theorem
\ref{evolve1} characterize the Cayley condition on
$\phi_1,\phi_2,\psi$ and that there is a natural conformal structure
on $\Sigma$.

\subsection{Gauge Transformations}
\label{gauge}

Let $\phi_1,\phi_2$ satisfy \eq{Cay7} and \eq{Cay8} in Theorem
\ref{evolve1} for some map $f$. Taking the triple cross product of
\eq{Cay7} and \eq{Cay8} with $\phi_1$ and $\phi_2$ gives:
\begin{align}
\label{Cay9} \frac{\partial\phi_1}{\partial s} & =
-\phi_1\times\phi_2\times\frac{\partial\phi_1}{\partial t}
+f^{'}\phi_2, \\
\label{Cay10} \frac{\partial\phi_2}{\partial s} & =
-\phi_1\times\phi_2\times \frac{\partial\phi_2}{\partial t}-
f^{'}\phi_1,
\end{align}
\noindent for some function $f^{'}:\Sigma\rightarrow\R$.

We are allowed to perform a rotation $\Theta(\sigma)$ to the
$(\phi_1(\sigma),\phi_2(\sigma))$-plane at each point
$\sigma\in\Sigma$ as long as the function $\Theta$ is smooth.  The
choice of $\Theta$ will then alter $f$ and $f^{'}$. We call such a
transformation a \emph{gauge transformation}.

We now show that under certain conditions there exists a gauge
transformation such that $f=f^{'}=0$.  Let
$\Theta:\Sigma\rightarrow\R$ be a smooth function and define
$\tilde{\phi}_1,\tilde{\phi}_2$ by
\begin{equation*}
\left(\begin{array}{c} \tilde{\phi}_1 \\ \tilde{\phi}_2
\end{array}\right)  =  \left(\begin{array}{rr} \cos\Theta &
\sin\Theta \\ -\sin\Theta & \cos\Theta\end{array}\right)
\left(\begin{array}{c} \phi_1 \\ \phi_2\end{array}\right).
\end{equation*}
\noindent Then $\tilde{\phi}_1,\tilde{\phi}_2$ satisfy \eq{Cay7} and
\eq{Cay8} with $f$ replaced by
$\tilde{f}=f+\frac{\partial\Theta}{\partial t}$.  Moreover, they
satisfy \eq{Cay9} and \eq{Cay10} with $f^{\prime}$ replaced by
$\tilde{f}^{\prime}=f^{\prime}+\frac{\partial\Theta}{\partial s}$.
Therefore, locally, there exists a smooth function $\Theta$ such
that $\tilde{f}=\tilde{f}^{\prime}=0$ if and only if $\frac{\partial
f}{\partial s}=\frac{\partial f^{\prime}}{\partial t}$.

If we differentiate \eq{Cay7} with respect to $s$ and
differentiate \eq{Cay9} with respect to $t$ we get:
\begin{align}
\label{double1} \frac{\partial^2\phi_1}{\partial s\partial t} & =
 \hspace{8pt}\phi_1\times\frac{\partial\phi_2}{\partial
s}\times\frac{\partial\phi_1}{\partial s} +
\phi_1\times\phi_2\times\frac{\partial^2\phi_1}{\partial s^2} +
\frac{\partial f}{\partial s}\hspace{3pt}\phi_2 +
f\hspace{2pt}\frac{\partial\phi_2}{\partial s}, \\
\label{double2} \frac{\partial^2\phi_1}{\partial t\partial s} & =
 -\phi_1\times\frac{\partial\phi_2}{\partial
t}\times\frac{\partial\phi_1}{\partial t} -
\phi_1\times\phi_2\times\frac{\partial^2\phi_1}{\partial t^2} +
\frac{\partial f^{\prime}}{\partial t}\phi_2 +
f^{\prime}\frac{\partial\phi_2}{\partial t}.
\end{align}
\noindent We must have that \eq{double1} and \eq{double2} are equal.
In particular, the inner products of $\phi_2$ with \eq{double1} and
\eq{double2} must be equal.  Note that
\begin{align*}
g\left(\phi_2,\frac{\partial^2\phi_1}{\partial s\partial t}\right)
& =  -\Phi\left(\phi_1,\phi_2,\frac{\partial\phi_1}{\partial
s},\frac{\partial\phi_2}{\partial s}\right) +\frac{\partial
f}{\partial s}, \\
g\left(\phi_2,\frac{\partial^2\phi_1}{\partial t\partial s}\right) &
=  \Phi\left(\phi_1,\phi_2,\frac{\partial\phi_1}{\partial
t},\frac{\partial\phi_2}{\partial t}\right) +\frac{\partial
f^{\prime}}{\partial t}
\end{align*}
\noindent and that
\begin{align}
\label{gauge1}
\Phi\left(\phi_1,\phi_2,\frac{\partial\phi_1}{\partial
s},\frac{\partial\phi_2}{\partial s}\right) & =
g\left(\frac{\partial\phi_1}{\partial
t},\frac{\partial\phi_2}{\partial s}\right), \\
\label{gauge2}
\Phi\left(\phi_1,\phi_2,\frac{\partial\phi_1}{\partial
t},\frac{\partial\phi_2}{\partial t}\right) & =
-g\left(\frac{\partial\phi_1}{\partial
s},\frac{\partial\phi_2}{\partial t}\right).
\end{align}
\noindent However,
\begin{equation*}
\Phi\left(\phi_1,\phi_2,\frac{\partial\phi_2}{\partial
s},\frac{\partial\phi_1}{\partial
s}\right)=g\left(\frac{\partial\phi_1}{\partial
s},\frac{\partial\phi_2}{\partial t}\right).
\end{equation*}
Hence, since $\Phi$ is alternating, we deduce that $\frac{\partial
f}{\partial s}=\frac{\partial f^{\prime}}{\partial t}$ if and only
if all the terms in \eq{gauge1} and \eq{gauge2} are zero.

We say that functions $\phi_1$ and $\phi_2$ satisfying \eq{Cay7},
\eq{Cay8}, \eq{Cay9} and \eq{Cay10} with $f=f^{\prime}=0$ are in the
\emph{flat gauge}.

We now give a geometric interpretation of the flat gauge.  Let
$(\Sigma,\pi)$ be a\linebreak 2-ruling.  Then, there is an
$\mathcal{S}^1$ bundle $\pi_P:P\rightarrow\Sigma$ as described after
Definition \ref{2Ruled}. An r-framing, which is equivalent to a
choice of $\phi_1,\phi_2$, gives a trivialization of $P$ and we can
consider it as a $\U(1)$ bundle.   Define a connection $\nabla_P$ on
$P$ by a connection 1-form given by $d\theta-f^{\prime}ds-fdt$,
where $\theta$ corresponds to the $\U(1)$ direction. This connection
is independent of the choice of r-framing, by the work above, and
has curvature which may be written as $(\frac{\partial
f^{\prime}}{\partial t}-\frac{\partial f}{\partial s})ds\w dt$.
Hence, the connection $\nabla_P$ defined by $\phi_1,\phi_2$ is flat
if and only if $\phi_1,\phi_2$ can be put in the flat gauge by some
gauge transformation locally.

\subsection{Planar 2-ruled Cayley 4-folds}

In this subsection we show that maps $\phi_1,\phi_2,\psi$ which do
not satisfy \eq{Cay7}-\eq{Caypsi} for any local oriented coordinates
$(s,t)$ on $\Sigma$ define a \emph{planar} Cayley 4-fold.

The next result shows that \eq{Cay7}-\eq{Caypsi} can be considered
as \emph{evolution equations} for $\phi_1,\phi_2,\psi$. We make the
definition here that a function is real analytic on a compact
interval $I$ in $\R$ if it extends to a real analytic function on an
open set containing $I$.

\begin{thm}
\label{evolve2} Let $I$ be a compact interval in $\R$, let $s$ be a
coordinate on $I$, let
$\phi_1^{\prime},\phi_2^{\prime}:I\rightarrow\mathcal{S}^7$ be
orthogonal real analytic maps and let
$\psi^{\prime}:I\rightarrow\R^8$ be a real analytic map. Let $N$ be
a neighbourhood of $0$ in $\R$ and let $f:I\times N\rightarrow\R$ be
a real analytic map.  Then there exist $\epsilon>0$ and unique real
analytic maps
$\phi_1,\phi_2:I\times(-\epsilon,\epsilon)\rightarrow\mathcal{S}^7$,
with $\phi_1,\phi_2$ orthogonal, and
$\psi:I\times(-\epsilon,\epsilon)\rightarrow\R^8$ satisfying
$\phi_1(s,0)=\phi_1^{\prime}(s)$, $\phi_2(s,0)=\phi_2^{\prime}(s)$,
$\psi(s,0)=\psi^{\prime}(s)$ for all $s\in I$ and
\begin{align}
\label{Cay7a}\frac{\partial\phi_1}{\partial t} & =
\phi_1\times\phi_2\times \frac{\partial\phi_1}{\partial
s} + f\phi_2, \\
\label{Cay8a} \frac{\partial\phi_2}{\partial t} & =
\phi_1\times\phi_2\times \frac{\partial\phi_2}{\partial s} - f
\phi_1, \\
\label{Cay9a} \frac{\partial\psi}{\partial t} & =
\phi_1\times\phi_2\times\frac{\partial\psi}{\partial s},
\end{align}
\noindent where $t$ is a coordinate on $(-\epsilon,\epsilon)$ and
the triple cross product is defined in \eq{Caytriplecross}. Let
$M$ be defined by:
\begin{equation}
\label{Cay10a} M =
\{r_1\phi_1(s,t)+r_2\phi_2(s,t)+\psi(s,t):(r_1,r_2)\in\R^2,\hspace{2pt}
s\in I,\hspace{2pt} t\in (-\epsilon,\epsilon)\}.
\end{equation}
\noindent Then $M$ is an r-framed 2-ruled Cayley 4-fold in $\R^8$.
\end{thm}

\begin{proof}
Since $I$ is compact and
$\phi_1^{\prime},\phi_2^{\prime},\psi^{\prime},f$ are real analytic,
we may apply the \emph{Cauchy--Kowalevsky Theorem} \cite[p.
234]{Racke} from the theory of partial differential \linebreak
equations to give unique functions $\phi_1,\phi_2,\psi:I\times
(-\epsilon,\epsilon)\rightarrow\R^8$ satisfying the initial
conditions and \eq{Cay7a}-\eq{Cay9a}. We must now show that
$\phi_1,\phi_2$ map to $\mathcal{S}^7$ and are orthogonal.

We first note that
\begin{align*}
\frac{\partial}{\partial t} \,g(\phi_1,\phi_2) & =
g\left(\frac{\partial\phi_1}{\partial
t},\phi_2\right)+g\left(\phi_1,\frac{\partial\phi_2}{\partial
t}\right) \\
& =  f\left(g(\phi_2,\phi_2)-g(\phi_1,\phi_1)\right), \\
\frac{\partial}{\partial t}\, g(\phi_1,\phi_1) & =
2fg(\phi_1,\phi_2), \\
\frac{\partial}{\partial t}\, g(\phi_2,\phi_2) & =
-2fg(\phi_1,\phi_2).
\end{align*}
\noindent Then $g(\phi_j,\phi_k)$ for $j,k=1,2$ are real analytic
functions satisfying this system of partial differential equations,
together with the initial conditions
$g(\phi_1,\phi_1)=g(\phi_2,\phi_2)=1$ and
$g(\phi_1,\phi_2)=\frac{\partial}{\partial t}g(\phi_j,\phi_k)=0$ at
$t=0$ given by assumption.  The functions
$g(\phi_1,\phi_1)=g(\phi_2,\phi_2)\equiv 1$ and
$g(\phi_1,\phi_2)\equiv 0$ also satisfy these equations and initial
conditions. It therefore follows from the Cauchy--Kowalevsky Theorem
that these two solutions must be locally equal and hence, for
$\epsilon>0$ sufficiently small, $|\phi_1|=|\phi_2|=1$ and
$\phi_1,\phi_2$ are orthogonal.

We conclude from Theorem \ref{evolve1} that $M$ is an r-framed
2-ruled Cayley 4-fold.

\end{proof}

Note that the Cayley 4-fold $M$ resulting from Theorem \ref{evolve2}
does not depend on the function $f$.

\medskip

Let $(\Sigma,\pi)$ and $(\tilde{\Sigma},\tilde{\pi})$ be 2-rulings
of a 4-fold in $\R^n$. We say that these 2-rulings are
\emph{distinct} if the families of affine 2-planes,
$\mathcal{F}_{\Sigma}=\{\pi^{-1}(\sigma):\sigma\in\Sigma\}$ and
$\mathcal{F}_{\tilde{\Sigma}}=\{\tilde{\pi}^{-1}(\tilde{\sigma}):
\tilde{\sigma}\in\tilde{\Sigma}\}$, are different.  If
$\mathcal{F}^n$ is the family of all affine 2-planes in $\R^n$ we
can consider $(\Sigma,\pi)$ as a map from $\Sigma$ to
$\mathcal{F}^n$ given by $\sigma\mapsto\pi^{-1}(\sigma)$ with image
$\mathcal{F}_{\Sigma}$.

Our next result is analogous to \cite[Theorem 6 part 2]{Bryant2},
which relates to ruled SL 3-folds.

\begin{prop}\label{1param}
A 2-ruled Cayley 4-fold in $\R^8$ which admits a real analytic one-parameter
 family of distinct real analytic 2-rulings is locally
isomorphic to an affine Cayley 4-plane in $\R^8$.
\end{prop}

\begin{proof}Let $\{(\Sigma_u,\pi_u):u\in\R\}$ be a real analytic
family of distinct real analytic\linebreak 2-rulings for a Cayley
4-fold $M$. Then there exists $p\in M$ such that
$\Pi_u=\pi_u^{-1}(\pi_u(p))$ is not constant as a 2-plane in $\R^8$.
Hence we have a one-parameter family of planes $\Pi_u\ni p$ in $M$
such that $\frac{d\Pi_u}{du}\neq 0$ for some $u$, i.e. such that
$\Pi_u$ changes nontrivially.  Therefore $\{\Pi_u:u\in\R\}$ is a
real analytic one-dimensional family of planes in $M$ containing
$p$. The total space of this family is a real analytic 3-fold $N$
contained in $M$.  Moreover, every plane in $M$ containing $p$ lies
in the affine Cayley 4-plane $p+T_pM$ and thus $N\subseteq p+T_pM$.
By Theorem \ref{realanal2}, $M$ and $p+T_pM$ must coincide on a
connected component of $M$. The result follows.
\end{proof}

We now give the result claimed at the start of the subsection, which
is analogous to the result \cite[Proposition 5.3]{Joyce2} for ruled
SL 3-folds.

\begin{prop}
\label{planar} Any r-framed 2-ruled Cayley 4-fold $(M,\Sigma,\pi)$
in $\R^8$ defined locally by maps $\phi_1,\phi_2,\psi$ which do not
satisfy \eq{Cay7}-\eq{Caypsi} for any local oriented coordinates
$(s,t)$ on $\Sigma$ is locally isomorphic to an affine Cayley
4-plane in $\R^8$.
\end{prop}

\begin{proof}
We may take the 2-ruling $(\Sigma,\pi)$ to be locally real analytic
since $M$ is real analytic by Theorem \ref{realanal1}. Let $I=[0,1]$
and let $\gamma:I\rightarrow\Sigma$ be a real analytic curve in
$\Sigma$. If we set $\phi_1^{\prime}(s)=\phi_1(\gamma(s))$,
$\phi_2^{\prime}(s)=\phi_2 (\gamma(s))$,
$\psi^{\prime}(s)=\psi(\gamma(s))$, then by Theorem \ref{evolve2} we
construct $\tilde{\phi}_1,\tilde{\phi}_2,\tilde{\psi}$ defining an
r-framed 2-ruled Cayley 4-fold $\tilde{M}$ satisfying
\eq{Cay7}-\eq{Caypsi} of Theorem \ref{evolve1}.  We have that
$M,\tilde{M}$ coincide in the real analytic 3-fold
$\pi^{-1}(\gamma(I))$, and hence, by Theorem \ref{realanal2}, they
must be locally equal.  Therefore, $M$ locally admits a 2-ruling
$(\tilde{\Sigma},\tilde{\pi})$ satisfying \eq{Cay7}-\eq{Caypsi} of
Theorem \ref{evolve1}, which must be distinct from $(\Sigma,\pi)$.

Using the notation introduced before Proposition \ref{1param}, the
families of affine\linebreak 2-planes $\mathcal{F}_{\Sigma}$ and
$\mathcal{F}_{\tilde{\Sigma}}$ coincide in the family of affine
2-planes defined by points on $\gamma$, denoted
$\mathcal{F}_{\gamma}$.  Using local real analyticity of the
families, either $\mathcal{F}_{\Sigma}$ is equal to
$\mathcal{F}_{\tilde{\Sigma}}$ locally or they only meet in
$\mathcal{F}_{\gamma}$ locally. The former possibility is excluded
because the 2-rulings $(\Sigma,\pi)$ and
$(\tilde{\Sigma},\tilde{\pi})$ are distinct and hence the latter is
true.

Let $\gamma_1$ and $\gamma_2$ be distinct real analytic curves near
$\gamma$ in $\Sigma$ defining 2-rulings $(\Sigma_1,\pi_1)$ and
$(\Sigma_2,\pi_2)$, respectively, as above.  Then
$\mathcal{F}_{\Sigma}\cap\mathcal{F}_{\Sigma_j}$ is locally equal to
$\mathcal{F}_{\gamma_j}$ for $j=1,2$.  Hence, the 2-rulings
$(\Sigma_1,\pi_1)$ and $(\Sigma_2,\pi_2)$ are not distinct (that is,
$\mathcal{F}_{\Sigma_1}=\mathcal{F}_{\Sigma_2}$) if and only if
$\mathcal{F}_{\gamma_1}=\mathcal{F}_{\gamma_2}$, which implies that
$\gamma_1=\gamma_2$.  Therefore, distinct curves near $\gamma$ in
$\Sigma$ produce different 2-rulings of $M$ and hence $M$ has
infinitely many 2-rulings.

Let $\{\gamma_u:u\in\R\}$ be a one parameter family of distinct real
analytic curves near $\gamma$ in $\Sigma$ with $\gamma_0=\gamma$.
Each curve in the family defines a distinct real analytic 2-ruling
$(\Sigma_u,\pi_u)$. Applying Proposition \ref{1param} gives the
result.
\end{proof}

Note that in the proof of Theorem \ref{evolve1} the condition
\eq{Caypsi} on $\psi$ was forced by the linear independence of the
derivatives of $\phi_1,\phi_2$.
 However, as we shall see in $\S$\ref{ruleds}, non-planar 2-ruled
4-folds can be constructed when the derivatives of $\phi_1,\phi_2$
are linearly dependent.

Proposition \ref{planar} tells us that for any \emph{non-planar}
2-ruled Cayley 4-fold $M$ defined by maps $\phi_1,\phi_2,\psi$ on
$\Sigma$ there exist locally oriented coordinates $(s,t)$ on
$\Sigma$ such that \eq{Cay7}-\eq{Caypsi} are satisfied. We shall see
in the next subsection that there is therefore a natural conformal
structure upon $\Sigma$, and $(s,t)$ are oriented conformal
coordinates with respect to this structure.

\subsection{Main Results}

We now present the main results of the paper on 2-ruled calibrated
4-folds in $\R^7$ and $\R^8$.  The first follows immediately from
Theorem \ref{evolve1} and Proposition \ref{planar}.

\begin{thm}
\label{Cayevolve} Let $(M,\Sigma,\pi)$ be a non-planar, r-framed,
2-ruled 4-fold in $\R^8$ defined by orthogonal real analytic maps
$\phi_1,\phi_2:\Sigma\rightarrow\mathcal{S}^7$ and a real analytic
map $\psi:\Sigma\rightarrow\R^8$ as follows:
\begin{equation}
\label{Cay1a} M =
\{r_1\phi_1(\sigma)+r_2\phi_2(\sigma)+\psi(\sigma):r_1,r_2\in\R,
\hspace{2pt}\sigma\in \Sigma\}.
\end{equation}
Then $M$ is Cayley if and only if there exist locally oriented
coordinates $(s,t)$ on $\Sigma$ such that
\begin{align}
\label{Cay2a} \frac{\partial\phi_1}{\partial t} & =
\phi_1\times\phi_2\times \frac{\partial\phi_1}{\partial
s} + f\phi_2, \\
\label{Cay3a} \frac{\partial\phi_2}{\partial t} & =
\phi_1\times\phi_2\times \frac{\partial\phi_2}{\partial s} - f
\phi_1, \\
\label{Cay4a} \frac{\partial\psi}{\partial t} & =
\phi_1\times\phi_2\times\frac{\partial\psi}{\partial
s}+g_1\phi_1+g_2\phi_2,
\end{align}
\noindent where the triple cross product is defined by equation
\eq{Caytriplecross} and $f,g_1,g_2:\Sigma\rightarrow\R$ are some
real analytic functions.
\end{thm}

We now prove the result claimed at the end of the last subsection.

\begin{prop} \label{conformal} Let $(M,\Sigma,\pi)$ be a non-planar, r-framed,
2-ruled Cayley\linebreak 4-fold in $\R^8$.  Then there exists a
unique conformal structure on $\Sigma$ with respect to which $(s,t)$
as given in Theorem \ref{Cayevolve} are oriented conformal
coordinates.
\end{prop}

\begin{proof}
Let $(s,t)$ be local oriented coordinates as given by Theorem
\ref{Cayevolve}.  Define a complex structure $J$ on $\Sigma$ by
requiring that $u=s+it$ is a holomorphic coordinate on $\Sigma$,
i.e. that $\frac{\partial}{\partial t}=J\frac{\partial}{\partial
s}$.  Note that $\phi_1,\phi_2$ as given in Theorem \ref{Cayevolve}
satisfy
\begin{align}
\label{Cay2b} \frac{\partial\phi_j}{\partial t}^{\perp} & =
\phi_1\times\phi_2\times \frac{\partial\phi_j}{\partial
s}^{\perp}, \\
\label{Cay3b} \frac{\partial\phi_j}{\partial s}^{\perp} & = -
\phi_1\times\phi_2\times \frac{\partial\phi_j}{\partial
t}^{\perp},
\end{align}
for $j=1,2$.  Suppose that $(\tilde{s},\tilde{t})$ are local
oriented coordinates on $\Sigma$ such that $\phi_1,\phi_2$ also
satisfy \eq{Cay2a}-\eq{Cay3a} in these coordinates.  Hence,
$\phi_1,\phi_2$ satisfy \eq{Cay2b}-\eq{Cay3b} for the coordinates
$(\tilde{s},\tilde{t})$.

We then calculate:
\begin{align*}
\frac{\partial\phi_j}{\partial\tilde{t}}^{\perp} & = \frac{\partial
s}{\partial\tilde{t}}\frac{\partial\phi_j}{\partial
s}^{\perp}+\frac{\partial
t}{\partial\tilde{t}}\frac{\partial\phi_j}{\partial t}^{\perp}
\\
&=\phi_1\times\phi_2\times\left(\frac{\partial
t}{\partial\tilde{t}} \frac{\partial\phi_j}{\partial
s}^{\perp}-\frac{\partial
s}{\partial\tilde{t}}\frac{\partial\phi_j}{\partial
t}^{\perp}\right) \\
\intertext{and} \frac{\partial\phi_j}{\partial\tilde{s}}^{\perp} &
= \frac{\partial
s}{\partial\tilde{s}}\frac{\partial\phi_j}{\partial
s}^{\perp}+\frac{\partial
t}{\partial\tilde{s}}\frac{\partial\phi_j}{\partial t}^{\perp}.
\end{align*}
Note that, from \eq{Cay2b}, $\frac{\partial\phi_j}{\partial
t}^\perp$ is orthogonal to $\frac{\partial\phi_j}{\partial s}^\perp$
and, moreover, that $\frac{\partial\phi_j}{\partial t}^\perp\neq0$
if and only if $\frac{\partial\phi_j}{\partial s}^\perp\neq0$ by the
definition of $f^\perp$ for a function $f:\Sigma\rightarrow\R$ and
the properties of the triple cross product. Using
\eq{Cay2b}-\eq{Cay3b} for $(\tilde{s},\tilde{t})$ we deduce that
\begin{equation}
\label{conf} \frac{\partial s}{\partial\tilde{s}}=\frac{\partial
t}{\partial\tilde{t}}\quad\text{and}\quad\frac{\partial
s}{\partial\tilde{t}}=-\frac{\partial t}{\partial\tilde{s}}\,,
\end{equation}
since not both $\frac{\partial\phi_1}{\partial
s}^{\perp},\frac{\partial\phi_2}{\partial s}^{\perp}$ are zero.

Therefore, using \eq{conf},
\begin{align*}
\frac{\partial}{\partial\tilde{t}}&=\frac{\partial
s}{\partial\tilde{t}}\frac{\partial}{\partial s}+\frac{\partial
t}{\partial\tilde{t}}\frac{\partial}{\partial t}=-\frac{\partial
t}{\partial\tilde{s}}\frac{\partial}{\partial s}+\frac{\partial
s}{\partial\tilde{s}}\frac{\partial}{\partial t}
=J\left(\frac{\partial
t}{\partial\tilde{s}}\frac{\partial}{\partial t}+\frac{\partial
s}{\partial\tilde{s}}\frac{\partial}{\partial s}\right)
=J\frac{\partial}{\partial\tilde{s}}.
\end{align*}
Hence we have the result.
\end{proof}

It is clear that the conformal structure given by Proposition
\ref{conformal} coincides with that given by the metric as
described in the preamble to Theorem \ref{evolve1}.

\medskip

We now use Proposition \ref{specialCay} in order to prove analogous
results for SL 4-folds in $\C^4$ and coassociative 4-folds in
$\R^7$.

We begin with the SL case and define the triple cross product of
$x,y,z$ in $\C^4$ by
\begin{equation}
\label{SLtriplecross} (x\times y\times z)^e =
(\Re\Omega)_{abcd}x^ay^bz^cg^{de},
\end{equation}
\noindent using index notation for tensors on $\C^4$, where $g^{de}$
is the inverse of the Euclidean metric on $\C^4$. By equation
\eq{CPhi}, this triple cross product agrees with the one defined in
\eq{Caytriplecross} when $\omega(x,y)=\omega(y,z)=\omega(z,x)=0$.

\begin{thm}
\label{SLevolve} Let $(M,\Sigma,\pi)$ be a non-planar, r-framed,
2-ruled 4-fold in $\C^4\cong\R^8$ defined by orthogonal real
analytic maps $\phi_1,\phi_2:\Sigma\rightarrow\mathcal{S}^7$ and a
real analytic map $\psi:\Sigma\rightarrow\R^8$ as follows:
\begin{equation}
\label{SL1} M =
\{r_1\phi_1(\sigma)+r_2\phi_2(\sigma)+\psi(\sigma):r_1,r_2\in\R,
\hspace{2pt}\sigma\in \Sigma\}.
\end{equation}
Then $M$ is special Lagrangian if and only if
$\omega(\phi_1,\phi_2)\equiv 0$ and there exist locally oriented
coordinates $(s,t)$ on $\Sigma$ such that:
\begin{align}
\label{SLcond1} \omega\left(\phi_j,\frac{\partial\phi_k}{\partial
s}\right)&\equiv 0& &\text{for $j,k=1,2$}, \\
\label{SLcond2} \omega\left(\phi_j,\frac{\partial\psi}{\partial
s}\right) &\equiv 0&
&\text{for $j=1,2$},\\
 \label{SL2} \frac{\partial\phi_1}{\partial t} & =
\phi_1\times\phi_2\times \frac{\partial\phi_1}{\partial
s} + f\phi_2,& &\\
\label{SL3} \frac{\partial\phi_2}{\partial t} & =
\phi_1\times\phi_2\times \frac{\partial\phi_2}{\partial s} - f
\phi_1,& &\\
\label{SL4} \frac{\partial\psi}{\partial t} & =
\phi_1\times\phi_2\times\frac{\partial\psi}{\partial
s}+g_1\phi_1+g_2\phi_2,&&
\end{align}
\noindent where the triple cross product is defined by equation
\eq{SLtriplecross}, and $f,g_1,g_2:\Sigma\rightarrow\R$ are some
real analytic functions.
\end{thm}

It is worth making clear that \eq{SL2}-\eq{SL4} are not the same
as \eq{Cay2a}-\eq{Cay4a} because of the different definitions of
the triple cross product.

\begin{proof}
By Proposition \ref{specialCay}, $M$ is SL if and only if $M$ is
Cayley and $\omega|_M\equiv 0$.  We thus conclude from Theorem
\ref{Cayevolve} that $M$ is SL if and only if $\phi_1,\phi_2,\psi$
satisfy \eq{Cay2a}-\eq{Cay4a} and $\omega|_{T_pM}\equiv 0$ for all
$p\in M$.  Therefore $\omega$ vanishes on $\langle
x,y,z,w\rangle_{\R}$, where $x=\phi_1(\sigma)$, $y=\phi_2(\sigma)$,
$z=r_1\frac{\partial\phi_1}{\partial s}(\sigma)+
r_2\frac{\partial\phi_2}{\partial
s}(\sigma)+\frac{\partial\psi}{\partial s} (\sigma)$,
$w=r_1\frac{\partial\phi_1}{\partial t}(\sigma)+
r_2\frac{\partial\phi_2}{\partial
t}(\sigma)+\frac{\partial\psi}{\partial t} (\sigma)$, for all
$(r_1,r_2)\in\R^2$, $\sigma\in\Sigma$.  Hence, the equations that
must be satisfied are $\omega(\phi_1,\phi_2)\equiv 0$ and
\begin{align}
\label{omega1}\omega\left(\phi_j, \frac{\partial\phi_k}{\partial
s}\right)=\omega\left(\phi_j,\frac{\partial\phi_k}{\partial
t}\right)&\equiv  0 & &\text{for $j,k=1,2$}, \\
\label{omega2} \omega\left(\phi_j,\frac{\partial\psi}{\partial
s}\right)=\omega\left(\phi_j,\frac{\partial\psi}{\partial
t}\right)&\equiv 0 & &\text{for $j=1,2$}, \\
\label{omega3}\omega\left(\frac{\partial\phi_j}{\partial
s},\frac{\partial\phi_j}{\partial
t}\right)=\omega\left(\frac{\partial\psi}{\partial
s},\frac{\partial\psi}{\partial t}\right)&\equiv  0 & &\text{for
$j=1,2$},
\\
\label{omega4}\omega\left(\frac{\partial\phi_1}{\partial
s},\frac{\partial\phi_2}{\partial
t}\right)+\omega\left(\frac{\partial\phi_2}{\partial
s},\frac{\partial\phi_1}{\partial t}\right)&\equiv  0, \\
\label{omega5}\omega\left(\frac{\partial\phi_j}{\partial
s},\frac{\partial\psi}{\partial
t}\right)+\omega\left(\frac{\partial\psi}{\partial
s},\frac{\partial\phi_j}{\partial t}\right)&\equiv  0 & &\text{for
$j=1,2$}.
\end{align}
However, if the functions $\phi_1,\phi_2,\psi$ satisfy
\eq{omega1}-\eq{omega5} and \eq{Cay2a}-\eq{Cay4a}, then they must
satisfy \eq{SL2}-\eq{SL4}.  Hence, it is enough to show that the
conditions in the theorem force \eq{omega1}-\eq{omega5} to hold in
order to prove the result.

If $x,y,z,w$ are vectors in $\C^4$ such that $\omega$ vanishes on
$\langle x,y,z,w \rangle_{\R}$, then direct calculation in
coordinates shows that
\begin{equation}
\label{perm} \omega(x, y\times z\times w) =
\Im(\epsilon_{abcd}x^ay^bz^cw^d),
\end{equation}
\noindent using index notation for tensors on $\C^4$, where
$\epsilon_{abcd}$ is the permutation symbol and the triple cross
product is given in \eq{SLtriplecross}. Noting that
$\omega(\phi_1,\phi_2)\equiv 0$, that \eq{SLcond1} and \eq{SLcond2}
hold, and the relationship between the triple cross products on
$\C^4$ and $\R^8$, we see that \eq{SL2}-\eq{SL4} hold.
 Hence $\omega(\phi_j,\frac{\partial\phi_k}{\partial t})=0$ for all
$j,k$ using \eq{SL2}-\eq{SL3} and \eq{perm}. Therefore \eq{omega1}
is satisfied.  Moreover, \eq{SL4} and \eq{perm} imply that
\eq{omega2} is satisfied.  If we use \eq{SLcond1}-\eq{SLcond2},
\eq{SL2}-\eq{SL4} and \eq{perm} again, we have that \eq{omega3} is
satisfied.

We now show that \eq{omega4} and \eq{omega5} are satisfied.
Calculation using \eq{SLcond1}, \eq{SL2}-\eq{SL3} and \eq{perm}
gives:
\begin{align*}
\omega\left(\frac{\partial\phi_1}{\partial
s},\frac{\partial\phi_2}{\partial
t}\right)&+\omega\left(\frac{\partial\phi_2}{\partial
s},\frac{\partial\phi_1}{\partial t}\right)
\\ &=\omega\left(\frac{\partial\phi_1}{\partial
s},\phi_1\times\phi_2\times \frac{\partial\phi_2}{\partial
s}\right)+\omega\left(\frac{\partial\phi_2}{\partial
s},\phi_1\times\phi_2\times \frac{\partial\phi_1}{\partial
s}\right) \\
&=\Im\left(\epsilon_{abcd}\frac{\partial\phi_1}{\partial
s}^a\phi_1^b\phi_2^c\frac{\partial\phi_2}{\partial
s}^d\right)+\Im\left(\epsilon_{abcd}\frac{\partial\phi_2}{\partial
s}^a\phi_1^b\phi_2^c\frac{\partial\phi_1}{\partial s}^d\right) \\
&=\Im\left((\epsilon_{abcd}+\epsilon_{dbca})\frac{\partial\phi_1}{\partial
s}^a\phi_1^b\phi_2^c\frac{\partial\phi_2}{\partial
s}^d\right)\equiv 0
\end{align*}
by the definition of the permutation symbol.  Hence \eq{omega4} is
satisfied.  An entirely similar argument using
\eq{SLcond1}-\eq{SL4} and \eq{perm} gives that \eq{omega5} is
satisfied.
\end{proof}

For the coassociative case we define the triple cross product of
$x,y,z$ in $\R^7$ by
\begin{equation}
\label{coasstriplecross} (x\times y\times
z)^e=(\ast\varphi)_{abcd}x^ay^bz^cg^{de}
\end{equation}
\noindent using index notation for tensors on $\R^7$, where $g^{de}$
is the inverse of the Euclidean metric on $\R^7$.  If we embed
$\R^7$ as $\{0\}\times\R^7$ in $\R^8$, then \eq{7Phi} implies that
this triple cross product agrees with \eq{Caytriplecross} when
$\varphi(x,y,z)=0$.

\begin{thm}
\label{coassevolve} Let $(M,\Sigma,\pi)$ be a non-planar,
r-framed, 2-ruled 4-fold in $\R^7$ defined by orthogonal real
analytic maps $\phi_1,\phi_2:\Sigma\rightarrow\mathcal{S}^6$, and
a real analytic map $\psi:\Sigma\rightarrow\R^7$ as follows:
\begin{equation}
\label{coass1} M =
\{r_1\phi_1(\sigma)+r_2\phi_2(\sigma)+\psi(\sigma):r_1,r_2\in\R,
\hspace{2pt}\sigma\in \Sigma\}.
\end{equation}
Then $M$ is coassociative if and only if there exist locally
oriented coordinates $(s,t)$ on $\Sigma$ such that
\begin{align}
\label{coasscond1}
\varphi\left(\phi_1,\phi_2,\frac{\partial\phi_j}{\partial
s}\right)&\equiv 0& &\text{for $j=1,2$,} \\
\label{coasscond2}\varphi\left(\phi_1,\phi_2,\frac{\partial\psi}{\partial
s}\right) &\equiv 0, && \\
 \label{coass2} \frac{\partial\phi_1}{\partial t} & =
\phi_1\times\phi_2\times \frac{\partial\phi_1}{\partial
s} + f\phi_2, &&\\
\label{coass3} \frac{\partial\phi_2}{\partial t} & =
\phi_1\times\phi_2\times \frac{\partial\phi_2}{\partial s} - f
\phi_1, &&\\
\label{coass4} \frac{\partial\psi}{\partial t} & =
\phi_1\times\phi_2\times\frac{\partial\psi}{\partial
s}+g_1\phi_1+g_2\phi_2,&&
\end{align}
\noindent where the triple cross product is defined by equation
\eq{coasstriplecross}, and $f,g_1,g_2:\Sigma\rightarrow\R$ are
some real analytic functions.
\end{thm}

\begin{proof} By Proposition \ref{specialCay}, $M$ is
coassociative if and only if $M\subseteq\R^7\subseteq\R^8$ is
Cayley.  We may deduce from Theorem \ref{Cayevolve} that $M$ is
Cayley if and only if there exist locally coordinates $(s,t)$ such
that $\phi_1,\phi_2,\psi$ satisfy \eq{Cay2a}-\eq{Cay4a}.  We then
note that \eq{coasscond1}-\eq{coasscond2} and the relationship
between the triple cross products on $\R^7$ and $\R^8$ ensure that
\eq{coass2}-\eq{coass4} are equivalent to \eq{Cay2a}-\eq{Cay4a}. The
result follows.
\end{proof}

\subsection{Holomorphic Vector Fields}

We now finish this section by giving a means of constructing
r-framed 2-ruled calibrated 4-folds starting from r-framed 2-ruled
calibrated cones using \emph{holomorphic vector fields}, which is
analogous to \cite[Theorem 6.1]{Joyce2}.

Suppose that $M_0$ is an r-framed, 2-ruled, Cayley cone in $\R^8$
defined by maps $\phi_1,\phi_2:\Sigma\rightarrow\R^8$ as in
\eq{2ruledcone}.  Then Proposition \ref{conformal} gives us a
conformal structure on $\Sigma$ related to $\phi_1,\phi_2$ and hence
we can consider $\Sigma$ as a \emph{Riemann surface}. Therefore
$\Sigma$ has a natural complex structure $J$ and we may define
oriented conformal coordinates $(s,t)$ on $\Sigma$. Suppose further
that $\phi_1,\phi_2$ are in the flat gauge. Hence the equations
$\phi_1,\phi_2$ satisfy are:
\begin{align}
\label{standard1} \frac{\partial\phi_j}{\partial t} & =
\phi_1\times\phi_2\times\frac{\partial\phi_j}{\partial s}
& &\text{for $j=1,2$,} \\
\label{standard2} \frac{\partial\phi_j}{\partial s} & =  -
\phi_1\times\phi_2\times\frac{\partial\phi_j}{\partial t} &
&\text{for $j=1,2$.}
\end{align}
We note from equations \eq{standard1} and \eq{standard2} that
there is a correspondence between ``$\phi_1\times\phi_2\times$"
and the complex structure $J$ on $\Sigma$.

\begin{thm} \label{Cayholo} Let $M_0$ be an r-framed, 2-ruled, Cayley
cone in $\R^8$ defined by maps
$\phi_1,\phi_2:\Sigma\rightarrow\mathcal{S}^7$ in the flat gauge,
where $\Sigma$ is a Riemann surface.  Let $w$ be a holomorphic
vector field on $\Sigma$ and define a map
$\psi:\Sigma\rightarrow\R^8$ by
$\psi=\mathcal{L}_{w}\phi_1+\mathcal{L}_{iw}\phi_2$, where
$\mathcal{L}_{w},\mathcal{L}_{iw}$ denote the Lie derivatives with
respect to $w,iw$.  Let $M$ be defined by $\phi_1,\phi_2,\psi$ as in
\eq{Cay1a}.  Then $M$ is an r-framed 2-ruled Cayley 4-fold in
$\R^8$.
\end{thm}

\begin{proof} We need to show that $\psi$ as defined satisfies \eq{Cay4a}.
If $w$ is identically zero, then $\psi$ trivially satisfies
\eq{Cay4a}. Therefore we need only consider the case where $w$ has
isolated zeros. Since the condition for $M$ to be Cayley is a closed
condition on $M$, it is sufficient to prove that \eq{Cay4a} holds at
any point $\sigma\in\Sigma$ with $w(\sigma)\neq 0$.

Let $\sigma\in\Sigma$ be such a point.  Then, since $w$ is a
holomorphic vector field, there exists an open set in $\Sigma$
containing $\sigma$ with oriented conformal coordinates $(s,t)$
such that $w=\frac{\partial}{\partial s}$, so that
$iw=\frac{\partial}{\partial t}$.  Hence
$\psi=\frac{\partial\phi_1}{\partial
s}+\frac{\partial\phi_2}{\partial t}$ in a neighbourhood of
$\sigma$.

Let $(e_1,\ldots,e_8)$ be an oriented orthonormal basis for $\R^8$
and $A=|\frac{\partial\phi_1}{\partial s}(\sigma)|$.  We transform
coordinates on $\R^8$ using Spin$(7)$ such that
$\phi_1(\sigma)=e_1$, $\phi_2(\sigma)=e_2$,
$\frac{\partial\phi_1}{\partial s}(\sigma)=Ae_3$ and
$\frac{\partial\phi_2}{\partial s}(\sigma)=a_1e_1+\ldots+a_8e_8$,
for some real constants $a_1,\ldots,a_8$.  Clearly, by
\eq{standard1}, $\frac{\partial\phi_1}{\partial t}(\sigma)=Ae_4$ and
hence $a_1=a_2=a_4=0$ by the orthogonality conditions imposed on
$\frac{\partial\phi_2}{\partial s}$ in the flat gauge.
Differentiating \eq{standard1} gives:
\begin{align*}
\frac{\partial^2\phi_1}{\partial s\partial t} &=
\phi_1\times\frac{\partial\phi_2}{\partial
s}\times\frac{\partial\phi_1}{\partial s} +
\phi_1\times\phi_2\times\frac{\partial^2\phi_1}{\partial s^2}\,, \\
\frac{\partial^2\phi_2}{\partial t^2} & =
\frac{\partial\phi_1}{\partial
t}\times\phi_2\times\frac{\partial\phi_2}{\partial
s}+\phi_1\times\frac{\partial\phi_2}{\partial
t}\times\frac{\partial\phi_2}{\partial
s}+\phi_1\times\phi_2\times\frac{\partial^2\phi_2}{\partial
t\partial s}\,.
\end{align*}
\noindent Therefore,
\begin{align*}
\frac{\partial\psi}{\partial t} & =  \frac{\partial}{\partial
t}\left(\frac{\partial\phi_1}{\partial
s}+\frac{\partial\phi_2}{\partial t}\right) \\
& =  \phi_1\times\frac{\partial\phi_2}{\partial
s}\times\frac{\partial\phi_1}{\partial s}
+\frac{\partial\phi_1}{\partial
t}\times\phi_2\times\frac{\partial\phi_2}{\partial
s}+\phi_1\times\frac{\partial\phi_2}{\partial
t}\times\frac{\partial\phi_2}{\partial s} \\
&+ \phi_1\times\phi_2\times\frac{\partial}{\partial
s}\left(\frac{\partial\phi_1}{\partial
s}+\frac{\partial\phi_2}{\partial t}\right).
\end{align*}
Calculation using \eq{Caytriplecross} and \eq{standard1} shows that
\begin{align*}
\phi_1\times\frac{\partial\phi_2}{\partial
s}\times\frac{\partial\phi_1}{\partial s}(\sigma) &=
A(a_7e_5-a_8e_6-a_5e_7+a_6e_8), \\
\frac{\partial\phi_1}{\partial
t}\times\phi_2\times\frac{\partial\phi_2}{\partial s}(\sigma) &=
A(-a_3e_1-a_7e_5+a_8e_6+a_5e_7-a_6e_8), \\
\phi_1\times\frac{\partial\phi_2}{\partial
t}\times\frac{\partial\phi_2}{\partial s}(\sigma) &=
-(a_3^2+a_5^2+a_6^2+a_7^2+a_8^2)e_2.
\end{align*}
\noindent We conclude that $\psi$ satisfies \eq{Cay4a} at $\sigma$
with $g_1(\sigma)=-Aa_3=-g(\frac{\partial\phi_1}{\partial
s}(\sigma),\frac{\partial\phi_2}{\partial s}(\sigma))$ and
$g_2(\sigma)=-|\frac{\partial\phi_2}{\partial s}(\sigma)|^2$.  By
the invariance of the metric and the triple cross product under
Spin$(7)$, and the discussion above, $\psi$ satisfies \eq{Cay4a} for
some $g_1,g_2:\Sigma\rightarrow\R$. Hence, using Theorem
\ref{Cayevolve}, the result follows.
\end{proof}

This result does not extend to the SL case in the way we might
expect. The construction starting with a 2-ruled SL cone $M_0$ will
generally produce a 2-ruled Cayley, but not SL, 4-fold $M$. The fact
that $M$ is Cayley follows trivially from Theorem \ref{Cayholo}, but
if we impose the condition $\omega|_M\equiv 0$, then $\phi_1,\phi_2$
must satisfy
\begin{align*}
\omega\left(\frac{\partial\phi_1}{\partial
s},\frac{\partial\phi_2}{\partial s}\right)
&=-\omega\left(\frac{\partial\phi_1}{\partial
t},\frac{\partial\phi_2}{\partial t}\right)=0, \\
\omega\left(\frac{\partial\phi_1}{\partial
s},\frac{\partial\phi_2}{\partial t}\right)&=
\hspace{8pt}\omega\left(\frac{\partial\phi_1}{\partial
t},\frac{\partial\phi_2}{\partial s}\right)=0,
\end{align*}
\noindent wherever $w\neq 0$.  At such a point, either all the
derivatives of $\phi_1,\phi_2$ are zero or at least one is non-zero.
In the first case both $\phi_1$ and $\phi_2$ are locally constant.
Otherwise, suppose without loss of generality that
$\frac{\partial\phi_1}{\partial s}\neq 0$ at a point $\sigma$ such
that $w(\sigma)\neq 0$. Then $\langle\frac{\partial\phi_1}{\partial
s},\frac{\partial\phi_1}{\partial t}\rangle_{\C}=\C^2$, and it is
orthogonal to $\langle\phi_1,\phi_2\rangle_{\C}=\C^2$ since $M$ is
SL and $\phi_1,\phi_2$ are in the flat gauge.  Therefore
$\frac{\partial\phi_2}{\partial
s}\in\langle\frac{\partial\phi_1}{\partial
s},\frac{\partial\phi_1}{\partial t}\rangle_{\C}$.  Note that
$g(\frac{\partial\phi_1}{\partial t},\frac{\partial\phi_2}{\partial
s})=\omega(\frac{\partial\phi_1}{\partial
t},\frac{\partial\phi_2}{\partial s})=0$ and
$\omega(\frac{\partial\phi_1}{\partial
s},\frac{\partial\phi_2}{\partial s})=0$. Hence there exists
$\theta\in\R$ such that $\cos\theta\frac{\partial\phi_1}{\partial
s}+\sin\theta\frac{\partial\phi_2}{\partial s}=0$.  Using
\eq{standard1}, $\cos\theta\frac{\partial\phi_1}{\partial
t}+\sin\theta\frac{\partial\phi_2}{\partial t}=0$.  Therefore,
$\cos\theta\phi_1+\sin\theta\phi_2$ is constant on a neighbourhood
of $\sigma$ and thus on the component of $\Sigma$ containing
$\sigma$. Hence we have the following result.

\begin{thm} \label{SLholo} Let $M_0$ be an r-framed, 2-ruled,
SL cone in $\C^4\cong\R^8$ defined by
$\phi_1,\phi_2:\Sigma\rightarrow\mathcal{S}^7$ in the flat gauge,
where $\Sigma$ is a Riemann surface.  Let $w$ be a holomorphic
vector field on $\Sigma$ and define $\psi:\Sigma\rightarrow\C^4$ by
$\psi=\mathcal{L}_{w}\phi_1+\mathcal{L}_{iw}\phi_2$, where
$\mathcal{L}_{w},\mathcal{L}_{iw}$ denote the Lie derivatives with
respect to $w,iw$.  Let $M$ be defined by $\phi_1,\phi_2,\psi$ as in
\eq{SL1}.  Then $M$ is an r-framed 2-ruled Cayley 4-fold in $\R^8$,
which is SL if and only if $w\equiv 0$ or there exists $\theta\in\R$
for each component $K$ of $\Sigma$ such that
$\cos\theta\phi_1+\sin\theta\phi_2$ is constant on $K$.
\end{thm}
We do, however, have a similar result to Theorem \ref{Cayholo} for
coassociative 4-folds.

\begin{thm} \label{coassholo} Let $M_0$ be an r-framed, 2-ruled,
coassociative cone in $\R^7$\linebreak defined by
$\phi_1,\phi_2:\Sigma\rightarrow\mathcal{S}^6$ in the flat gauge,
where $\Sigma$ is a Riemann surface.  Let $w$ be a holomorphic
vector field on $\Sigma$ and define $\psi:\Sigma\rightarrow\R^7$ by
$\psi=\mathcal{L}_{w}\phi_1+\mathcal{L}_{iw}\phi_2$, where
$\mathcal{L}_{w},\mathcal{L}_{iw}$ denote the Lie derivatives with
respect to $w,iw$.  Let $M$ be defined by $\phi_1,\phi_2,\psi$ as in
\eq{coass1}.  Then $M$ is an r-framed 2-ruled coassociative 4-fold
in $\R^7$.
\end{thm}

\begin{proof} This follows immediately from Theorem
\ref{Cayholo} since $M\subseteq\R^7\subseteq\R^8$ is Cayley and
therefore coassociative by Proposition \ref{specialCay}.
\end{proof}

\section{Examples}
\label{exs}

We shall now exhibit explicit examples of 2-ruled 4-folds.
\subsection{$\U(1)$-Invariant 2-ruled Cayley 4-folds}

We consider the family of SL 4-folds in $\C^4$ given in
\cite[Theorem III.3.1]{HarLaw}.  Let ${\bf
c}=(c_1,c_2,c_3,c_4)\in\R^4$ be constant and define $M_{\bf
c}\subseteq\C^4$ by:
\begin{equation}
\label{SLHarLaw} M_{\bf c} =
\{(z_1,z_2,z_3,z_4)\in\C^4:\Re(z_1z_2z_3z_4)=c_1,\hspace{2pt}
|z_1|^2-|z_j|^2=c_j\; {\rm for} \, j=2,3,4\}.
\end{equation}
\noindent  Then $M_{\bf c}$ is an SL 4-fold in $\C^4$ invariant
under $\U(1)^3$.

Taking ${\bf c}=0$, we see that $M_0$ is an r-framed 2-ruled SL cone
in $\C^4$, with three different 2-rulings. For each of the distinct
2-rulings we are then able to apply the holomorphic vector field
result of Theorem \ref{SLholo} to obtain families of r-framed
2-ruled Cayley 4-folds which are invariant under $\U(1)$.

\begin{thm}
\label{example1} Let $w:\C\rightarrow\C$ be a holomorphic
function. Then
\begin{align}
M_1=
\bigg\{&\frac{1}{2}\Big(ie^{is}\left(re^{i\theta}\!+i\bar{w}(s+it)\right),
e^{-is}\left(re^{i\theta}\!-i\bar{w}(s+it)\right), \nonumber\\
\label{exCay1} & e^{it}\left(re^{-i\theta}\!+w(s+it)\right),
e^{-it}\left(re^{-i\theta}\!-w(s+it)\right)\Big):r,s,t,\theta\in\R\bigg\},
\\
M_2=
\bigg\{&\frac{1}{2}\Big(ie^{is}\left(re^{i\theta}\!+i\bar{w}(s+it)\right),
e^{-it}\left(re^{-i\theta}\!-w(s+it)\right), \nonumber\\
\label{exCay2}& e^{it}\left(re^{-i\theta}\!+w(s+it)\right),
e^{-is}\left(re^{i\theta}\!-i\bar{w}(s+it)\right)\Big):r,s,t,\theta\in\R\bigg\},
\\
M_3=
\bigg\{&\frac{1}{2}\Big(ie^{is}\left(re^{i\theta}\!+i\bar{w}(s+it)\right),
e^{it}\left(re^{-i\theta}\!+w(s+it)\right), \nonumber\\
\label{exCay3} &
e^{-is}\left(re^{i\theta}\!-i\bar{w}(s+it)\right),
e^{-it}\left(re^{-i\theta}\!-w(s+it)\right)\Big):r,s,t,\theta\in\R\bigg\}
\end{align}
\noindent are r-framed 2-ruled Cayley 4-folds in $\R^8\cong\C^4$.
\end{thm}

\begin{proof}
We only prove the result for $M_1$ as the proof for the other two is
similar.  In this example we define $M_0$ by functions
$\phi_1,\phi_2:\R^2\rightarrow\mathcal{S}^7\subseteq\C^4$ given by:
\begin{align}
\label{exCay1phi1}
\phi_1(s,t)&=\frac{1}{2}(ie^{is},e^{-is},\hspace{2pt}e^{it},\hspace{2pt}e^{-it}), \\
\label{exCay1phi2}
\phi_2(s,t)&=\frac{i}{2}(ie^{is},e^{-is},-e^{it},-e^{-it}),
\end{align}
\noindent so that $M_0$ is 2-ruled by planes of the form:
\begin{align*}
\Pi_{r,\theta} & =
\{r\cos\theta\phi_1(s,t)+r\sin\theta\phi_2(s,t):s,t\in\R\} \\
& =
\left\{\frac{r}{2}\left(ie^{i(\theta+s)},e^{i(\theta-s)},e^{-i(\theta-t)},e^{-i(\theta+t)}\right):s,t\in\R\right\}.
\end{align*}
We verify through direct calculation that
$g(\phi_1,\phi_2)=\omega(\phi_1,\phi_2)=0$, \eq{SLcond1}, \eq{SL2}
and \eq{SL3} are satisfied for $f=0$, and that $\phi_1,\phi_2$ are
in the flat gauge.

Let $w(s+it)=u(s,t)+iv(s,t)$ for functions $u,v:\R^2\rightarrow\R$.
Then $\psi=\mathcal{L}_{w}\phi_1+\mathcal{L}_{iw}\phi_2$ is given
by:
\begin{align}
&\psi(s,t)=u(s,t)\frac{\partial\phi_1}{\partial
s}+v(s,t)\frac{\partial\phi_1}{\partial
t}-v(s,t)\frac{\partial\phi_2}{\partial
s}+u(s,t)\frac{\partial\phi_2}{\partial t} \nonumber\\
\label{exCay1psi} &=
\frac{1}{2}\left(i\bar{w}(s+it)ie^{is},-i\bar{w}(s+it)e^{-is},w(s+it)e^{it},-w(s+it)e^{-it}\right).
\end{align}
\noindent We check that $\psi$ satisfies \eq{Cay4a} of Theorem
\ref{Cayevolve} with $g_1=0$ and $g_2=-\frac{1}{2}$, which agrees
with the calculations in the proof of Theorem \ref{Cayholo}.
Therefore, by Theorem \ref{Cayholo}, $M_1$ as defined in \eq{exCay1}
is an r-framed 2-ruled Cayley 4-fold.
\end{proof}

Note, from the proof above, that
$2\omega(\phi_1,\frac{\partial\psi}{\partial s})=v$ and
$2\omega(\phi_1,\frac{\partial\psi}{\partial s})=u$, so that
\eq{SLcond2} of Theorem \ref{SLevolve} is satisfied if and only if
$w\equiv 0$.  Therefore, if $w$ is not identically zero, Theorem
\ref{SLholo} shows that $M_1$ is an r-framed 2-ruled Cayley 4-fold
which is not special Lagrangian. Similarly for $M_2$ and $M_3$.

An interesting special case is when $w$ in Theorem \ref{example1} is
taken to be constant.  Here, calculation shows that, each $M_j$ is
invariant under a $\U(1)^2$ subgroup of $\U(1)^3$. Moreover, they
are asymptotically conical to $M_0$ with order $O(r^{-1})$, in the
sense of Definition \ref{asym}.

\subsection{Ruled Associative and Special Lagrangian 3-folds}
\label{ruleds}

We can construct examples of 2-ruled 4-folds from \emph{ruled}
associative 3-folds in $\R^7$ and SL 3-folds in $\C^3$, as described
in \cite[$\S6$]{Lotay} and \cite{Joyce2}.  We first give the
definitions, of which Definition \ref{2Ruled} is the analogue for
2-ruled 4-folds.

\begin{dfn}
\label{ruled} Let $M$ be a 3-dimensional submanifold of $\R^7$ (or
$\C^3$). A \emph{ruling} of $M$ is a pair $(\Sigma,\pi)$, where
$\Sigma$ is a 2-dimensional manifold and $\pi:M \rightarrow\Sigma$
is a smooth map, such that for all $\sigma\in\Sigma$ there exist
${\bf v}_{\sigma}\in \mathcal{S}^6$ (or $\mathcal{S}^5$), ${\bf
w}_{\sigma}\in\R^7$ (or $\C^3$) such that $\pi^{-1}(\sigma) =\{r{\bf
v}_{\sigma}+{\bf w}_{\sigma}:r\in\R\}$. Then the triple $(M,\Sigma,
\pi)$ is a \emph{ruled 3-fold} in $\R^7$ (or $\C^3$).

An \emph{r-orientation} for a ruling $(\Sigma,\pi)$ is a choice of
orientation for the affine straight line $\pi^{-1}(\sigma)$ in
$\R^7$ (or $\C^3$), for each $\sigma\in\Sigma$, which varies
continuously with $\sigma$.  Then a ruled 3-fold $(M,\Sigma,\pi)$
with an r-orientation is called an \emph{r-oriented ruled 3-fold}.

Let $(M,\Sigma,\pi)$ be an r-oriented ruled 3-fold.  For each
$\sigma \in\Sigma$ define $\phi(\sigma)$ to be the unique unit
vector in $\mathcal{S}^6$ (or $\mathcal{S}^5$) parallel to
$\pi^{-1}(\sigma)$ and in the positive direction with respect to the
orientation on $\pi^{-1}(\sigma)$, given by the r-orientation. Then
$\phi:\Sigma\rightarrow \mathcal{S}^6$ (or $\mathcal{S}^5$) is a
smooth map. Define $\psi:\Sigma \rightarrow\R^7$ (or $\C^3$) such
that, for all $\sigma\in\Sigma$, $\psi(\sigma)$ is the unique vector
in $\pi^{-1}(\sigma)$ orthogonal to $\phi(\sigma)$.  Then $\psi$ is
smooth and
\begin{equation}
\begin{split}
\label{ruled1}
M=\{r\phi(\sigma)+\psi(\sigma):\sigma\in\Sigma,\hspace{2pt}
r\in\R\}.
\end{split}
\end{equation}
\end{dfn}

\vspace{-10pt}

In the ruled case, there is a natural way to define a metric on
$\Sigma$ as the pullback $\phi^{\ast}(g)$ of the round metric $g$
on $\mathcal{S}^6$ (or $\mathcal{S}^5$), and hence we may always
define oriented conformal coordinates in terms of a natural
complex structure on $\Sigma$.

Suppose that $(N,\Sigma,\pi)$ is a ruled 3-fold.  Let $M=\R\times N$
and let $\tilde{\pi}:M\rightarrow\Sigma$ be given by
$\tilde{\pi}(r,p)=\pi(p)$ for all $p\in N$.  Clearly,
$(M,\Sigma,\tilde{\pi})$ is a 2-ruled 4-fold since
$\tilde{\pi}^{-1}(\sigma)=\R\times\pi^{-1}(\sigma)$ for all
$\sigma\in\Sigma$.  Suppose further that $(N,\Sigma,\pi)$ is
r-oriented.  Using the r-orientation, we have a natural choice of
oriented orthonormal basis for the plane $\tilde{\pi}^{-1}(\sigma)$,
which varies smoothly with $\sigma$.  Therefore,
$(M,\Sigma,\tilde{\pi})$ is r-framed.

We now state and prove the following theorem.

\begin{thm}\begin{itemize}
\item[\emph{(a)}] Let $N\subseteq\R^7$ be an (r-oriented) ruled associative 3-fold.
Then $\R\times N\subseteq\R\oplus\R^7\cong\R^8$ is an (r-framed)
2-ruled Cayley 4-fold.
\item[\emph{(b)}] Let $L\subseteq\C^3$ be an (r-oriented) SL 3-fold with phase
$-i$.  Then $\R\times L\subseteq\R\oplus\C^3\cong\R^7$ is an
(r-framed) 2-ruled coassociative 4-fold.
\end{itemize}
\end{thm}

\begin{proof}
Let $N$ be an associative 3-fold in $\R^7$.  Then, by equation
\eq{7Phi}, since $\ast\varphi|_{\R\times N}\equiv0$ and $N$ is a
$\varphi$-submanifold, it is clear that $\R\times N$ is calibrated
with respect to $\Phi$.  The comments before the theorem then give
the result (a).

Let $L$ be an SL 3-fold with phase $-i$ in $\C^3$.  Therefore, $L$
is calibrated with respect to $-\Im\Omega$, where $\Omega$ is the
holomorphic volume form on $\C^3$.  Let $(z_1,z_2,z_3)$ be complex
coordinates on $\C^3$.  We identify $\R\oplus\C^3$ and $\R^7$ by
defining coordinates on $\R^7$ as $(x_1,\ldots,x_7)$ where $x_1$ is
the coordinate on $\R$ and we let $z_1=x_2+ix_3$, $z_2=x_4+ix_5$,
$z_3=x_6+ix_7$.  Then,
\begin{equation}\label{Cstarphi}
*\varphi=\frac{1}{2}\,\omega\w\omega-dx_1\w\Im\Omega.
\end{equation}
Clearly, $\omega\w\omega|_{\R\times L}\equiv0$ as $\omega|_L\equiv0$
and hence $\R\times L$ is coassociative.  Again, we use the results
before the theorem to give the result (b).
\end{proof}

\subsection{Complex Cones}

We define a complex cone $C$ in
$\C^4$ by
\begin{equation*}
C=\{(z_1,z_2,z_3,z_4)\in\C^4:P(z_1,z_2,z_3,z_4)=Q(z_1,z_2,z_3,z_4)=0\}
\end{equation*}
\noindent where $P,Q$ are homogeneous complex polynomials.  Suppose
further that $P,Q$ are such that $C$ is non-planar and nonsingular
except at $0$. Define a projection $\tilde{\pi}$ from $C\setminus 0$
to $\C\P^3$ by $\tilde{\pi}((z_1,z_2,z_3,z_4))=[z_1,z_2,z_3,z_4]$
and let $\Sigma$ be the image of $\tilde{\pi}$.  Let $M_0$ be given
by
\begin{equation*}
M_0=\{(z_1,z_2,z_3,z_4,\sigma)\in
C\times\Sigma:(z_1,z_2,z_3,z_4)\in\sigma\}
\end{equation*}
and define $\iota:M_0\rightarrow\C^4$ by
$\iota(z_1,z_2,z_3,z_4,\sigma)=(z_1,z_2,z_3,z_4)$.  Then $\iota$ is
an immersion except at $0$ and thus $M_0$ can be considered as an
immersed submanifold of $\C^4$ which is only singular at $0$. Let
$\pi:M_0\rightarrow\Sigma$ be given by
$\pi(z_1,z_2,z_3,z_4,\sigma)=\sigma$. Clearly, $M_0$ is 2-ruled by
complex lines $\pi^{-1}(\sigma)$ in $\C^4$.  Since any complex
surface in $\C^4\cong\R^8$ is Cayley by \cite[$\S$IV.2.C]{HarLaw},
$(M_0,\Sigma,\pi)$ is a 2-ruled Cayley 4-fold.

We can define a local holomorphic coordinate $w$, and hence oriented
conformal coordinates $(s,t)$, in $\Sigma$ by $w\mapsto
[z_1,z_2,z_3,z_4](w)$ in some open set $U$ in $\Sigma$. Suppose
without loss of generality that $z_4\neq 0$ in $U$.  Then, we may
rescale so that $z_4=1$ and define maps
$\phi_1,\phi_2:U\rightarrow\mathcal{S}^7$ by
\begin{align*}
\phi_1(s,t)&=\left(\frac{z_1(s,t)}{r}\,, \frac{z_2(s,t)}{r}\,,
\frac{z_3(s,t)}{r}\,, \frac{1}{r}\right),
\\
\phi_2(s,t)&=i\phi_1(s,t),
\end{align*}
\noindent where
$r=(1+|z_1(s,t)|^2+|z_2(s,t)|^2+|z_3(s,t)|^2)^{\frac{1}{2}}$.  We
can then write $M_0$ locally in the form \eq{2ruledcone}, where
$\phi_1,\phi_2$ satisfy \eq{Cay2a}-\eq{Cay3a}, since $C$, and hence
$M_0$, is non-planar. If we define $M$ by \eq{Cay1a} where $\psi$
satisfies \eq{Cay4a} then, from Theorem \ref{Cayevolve}, $M$ is a
non-planar, r-framed, 2-ruled Cayley 4-fold in $\R^8$.

\appendix
\section{Appendix}

\subsection{Cayley Multiplication Table for the Octonions}
\label{octotable}

Let $e_1=1$ and let $\{e_2,\ldots,e_8\}$ be a basis for Im $\O$.
Then a Cayley multiplication table for the octonions is as shown:

\medskip

\begin{center}
$\begin{array}{rrrrrrrrrr}
  & \vline & e_1  &  e_2 &  e_3 &  e_4 &  e_5 &  e_6 &  e_7 &  e_8 \\
\hline
  e_1 & \vline &  e_1  &  e_2 &  e_3 &  e_4 &  e_5 &  e_6 &  e_7 &  e_8 \\
e_2 & \vline & e_2 &  -e_1  &  e_4 & -e_3 &  e_6 & -e_5 &  e_8 & -e_7 \\
e_3 & \vline & e_3 & -e_4 &  -e_1  &  e_2 &  e_7 & -e_8 & -e_5 &  e_6 \\
e_4 & \vline & e_4 &  e_3 & -e_2 &  -e_1  & -e_8 & -e_7 &  e_6 &  e_5 \\
e_5 & \vline & e_5 & -e_6 & -e_7 &  e_8 &  -e_1  &  e_2 &  e_3 & -e_4 \\
e_6 & \vline & e_6 &  e_5 &  e_8 &  e_7 & -e_2 &  -e_1  & -e_4 & -e_3 \\
e_7 & \vline & e_7 & -e_8 &  e_5 & -e_6 & -e_3 &  e_4 &  -e_1  &  e_2 \\
e_8 & \vline & e_8 &  e_7 & -e_6 & -e_5 &  e_4 &  e_3 & -e_2 &
-e_1
\end{array}$
\end{center}

\medskip

\noindent Note that the multiplication is defined above so as to be
compatible with the formula \eq{phi} for $\varphi$ and hence the
table is not the standard one.

\subsection{Calculating the Fourfold Cross Product}
\label{calc}

At various stages in the proof of Theorem \ref{evolve1} we are
required to make calculations involving the fourfold cross product
on $\O\cong\R^8$.  We shall give here details of these
calculations and the methods employed in order to compute these
products efficiently.

For our purposes we need only consider fourfold cross products of
the form
$$f_{jk}=e_1\times e_2\times e_j\times e_k.$$
It is clear by Definition \ref{defcross} that $f_{jk}$ is
antisymmetric and that $f_{jk}=0$ for $1\leq j,k\leq 2$ since the
fourfold cross product is alternating.

We only want to consider the case when $\text{Im }f_{jk}\neq 0$.
By Proposition \ref{Cay4plane} this occurs if and only if
$\{e_1,e_2,e_j,e_k\}$ does not lie in a Cayley 4-plane.  Hence we
deduce that $\text{Im }f_{jk}= 0$ if $\{j,k\}=\{3,4\}$ or
$\{j,k\}=\{5,6\}$ or $\{j,k\}=\{7,8\}$.

We next make the following observation.  By the invariance of the
fourfold cross product under Spin$(7)$, if $\{e_j,e_k,e_l,e_m\}$
is an ordered basis for a Cayley 4-plane, then either
$f_{jk}=f_{lm}$ or $f_{jk}=-f_{lm}$, depending on whether
$\{e_j,e_k,e_l,e_m\}$ is a positively oriented basis or not.

Therefore, the only fourfold cross products we require are:
\begin{align*}
f_{58}&= e_3=f_{67}, \\
f_{57}&= e_4=f_{86}, \\
f_{74}&= e_5=f_{83}, \\
f_{48}&= e_6=f_{73}, \\
f_{36}&= e_7=f_{45}, \\
f_{35}&= e_8=f_{64}.
\end{align*}




\end{document}